\documentclass[10pt,fleqn,leqno]{article}
\usepackage{amsbsy}
\usepackage{amsfonts}
\usepackage{latexsym}
\usepackage{graphicx}

\newcommand{\be}{\begin{equation}}
\newcommand{\ee}{\end{equation}}
\newcommand{\ba}{\begin{array}}
\newcommand{\ea}{\end{array}}

\newcommand{\sn}{{\rm sn}}
\newcommand{\cn}{{\rm cn}}
\newcommand{\dn}{{\rm dn}}
\newcommand{\kp}{{\mathfrak{K}}}

\newcommand{\bt}{\beta}

\newcommand{\bea}{\begin{eqnarray}}
\newcommand{\eea}{\end{eqnarray}}

\newcommand{\ga}{\gamma}

\newcommand{\La}{\Lambda}
\newcommand{\al}{\alpha}
\newcommand{\pri}{\prime}

\newcommand{\ce}{{\cal E}}
\newcommand{\cet}{\tilde{\cal E}}
\newcommand{\vt}{\vartheta}

\newcommand{\up}{{\rm {\bf u}}^{+}}
\newcommand{\um}{{\rm {\bf u}}^{-}}
\newcommand{\ux}{{\rm {\bf u}}_{x}}
\newcommand{\nbh}{n\hat{{\bf B}}}
\newcommand{\bh}{\hat{{\bf B}}}
\newcommand{\pg}{\mathfrak{p}}
\newcommand{\BZ}{\mathbb{Z}}

\newcommand{\wh}{\widehat}
\begin{document}
\newtheorem{theorem}{Theorem}[section]
\newtheorem{pro}[thm]{Proposition}
\newtheorem{lem}[thm]{Lemma}
\newtheorem{cor}[thm]{Corollary}

\title{Generalizations of  Chebyshev polynomials and  Polynomial 
Mappings}

\author{Yang Chen$^{\dag}$ and James Griffin$^{\dag\dag}$\\
Department of Mathematics, Imperial College,\\
180 Queen's Gate, London, SW7 2BZ, UK \\
Mourad E H Ismail$^{\dag\dag\dag}$\\
Department of Mathematics, University of Central Florida,\\
Orlando, Florida, 32816, USA}
\date{}

\maketitle
 
\begin{abstract} 

In this paper we show how polynomial mappings of degree 
$\kp$ from a union of 
disjoint intervals onto $[-1,1]$ generate a countable 
number of special 
cases of   generalizations of  Chebyshev polynomials. 
We also derive a new 
expression for these generalized  Chebyshev polynomials 
for any genus $g$, from which 
the coefficients of $x^n$ can be found explicitly in 
terms of the branch points 
and the recurrence coefficients. We find that this 
representation is useful for
 specializing to polynomial mapping cases for small 
$\kp$ where we will have explicit 
expressions for the recurrence coefficients in terms of 
the branch points. We 
study in detail certain special cases of the polynomials 
for small degree 
mappings and prove a theorem concerning the location of 
the zeroes of the 
polynomials. We also derive an explicit expression for 
the discriminant for 
the genus 1 case of our  Chebyshev polynomials that is 
valid for any configuration of the 
branch point.

\end{abstract}

\bigskip
{\bf Running title}: Chebyshev Polynomials, Generalized  .
\vskip .3cm $\dag\;$e-mail: y.chen@ic.ac.uk
\vskip .2cm $\dag\dag\;$e-mail:j.c.griffin@ic.ac.uk
\vskip .2cm $\dag\dag\dag\;$e-mail :ismail@math.usf.edu

\bigskip

\setcounter{equation}{0}
\setcounter{thm}{0}

\section{Introduction and Preliminaries.} 

Akhiezer \cite{Akh1}, \cite{Akh2} and, Akhiezer and 
Tom\v cuk \cite{Akh-Tom}  introduced orthogonal polynomials on two intervals
which generalize the Chebyshev polynomials. He observed that
the study of properties of these polynomials requires the use
of elliptic functions. In the case of more than two intervals, Tom\v cuk \cite{Tom}, investigated
their Bernstein-Szeg\H{o} asymptotics, with the theory of Hyperelliptic
integrals, and found expressions in terms of a certain Abelian integral of the third kind. 
However, in his formulation, certain unknown points on a Hyperelliptic Riemann surface emerge
due to the lack of an explicit representation of the original polynomials. This was
circumvented in \cite{Chen-Lawrence}. In his book on elliptic functions \cite{Akhbook}
Akhiezer, obtained explicit formulas for the two interval case as an example of the 
application of elliptic functions to the theory of conformal mapping.  

In 1984 Al-Salam, Allaway, and Askey
introduced sieved ultraspherical polynomials which are
orthogonal with respect to an absolutely continuous measure supported on 
$[-1,1]$ but the weight function vanishes at $k+1$ points. 
Ismail \cite{Ism2} 
observed that the vanishing of the weight function means that 
the polynomials  are orthogonal    
on several adjacent intervals. He then  introduced
one  additional parameter in the definition of the sieved
ultraspherical polynomials which made them orthogonal
on several intervals with gaps. In particular his
polynomials  include analogues of the Chebyshev polynomials
of the first and second kind.  Their continuous spectrum is
$T_k^{-1}([-c, c])$, for  $c \in (-1, 1)$, where $T_k$ is a
Chebyshev polynomial of the first kind. More over these polynomials have simple closed form 
expressions 
and elementary generating functions.  More general sieved polynomials are in \cite{Cha:Ism3}. 

In \cite{Peher1} a study of Chebyshev type polynomials was carried out. They arose as  
extremal polynomials and are orthogonal  on several intervals.  
The Chebyshev type polynomials satisfy similar extremal 
properties 
to the classical Chebyshev polynomials on $[-1,1]$. 
In particular they deviate least from zero on a union 
of disjoint intervals. The extremal polynomials also have 
the property that they are 
orthogonal with respect to some weight function of the form,
\bea
\frac{\pi_{g-1}(x)}{\sqrt{\prod_{j=1}^{2g}(x-a_j)}}, 
\quad x \in E, \nonumber
\eea
where $E := \cup_{j=1}^{2g-1}(a_j,a_{j+1})$ and the 
polynomial $\pi_{g-1}(x)$ is chosen 
such that the weight is positive on the interior of $E$. 
The existence of an extremal polynomials 
on any set $E$ is clearly not guaranteed. It was shown 
in \cite{Peher1} that the existence of an 
extremal polynomial on a set $E$ is equivalent to the 
asymptotic periodicity of the recurrence 
coefficients of a sequence of orthogonal polynomials on $E$. 
The extremal polynomial also has the
property that it maps the set $E$ onto $[-1,1]$.  

In \cite{Ger:Van}, using a polynomial mapping of this 
kind, a prescription was given for the 
construction of a sequence of orthogonal polynomials on 
$E$ from a sequence orthogonal on a single interval. 
The corresponding weight supported on $E$ 
can be obtained using the following procedure, 
outlined in \cite{Ger:Van}.  In this generality closed 
form expressions and generating functions are not available. 
This raises the question of finding a class of polynomials as general 
as possible which has the extremal properties of the polynomials  
in Peherstorfer \cite{Peher1} and the simple structural formulas of 
Ismail's polynomials in \cite{Ism2}.  The first step was carried out in 
 \cite{Chen-Lawrence} and is further developed here. The present paper we 
develop this  work further by borrowing ideas from polynomial mappings and 
sieved polynomials.  

We start by choosing a polynomial of degree $\kp-1$ 
and call this polynomial $W_{\kp-1}(t)$. Then we 
form the expression,
\bea
\frac{W_{\kp-1}(t)}{M_\kp(t)-x}=\frac{W_{\kp-1}(t)}
{l_{\kp}\prod_{i=1}^\kp(t-M_i^{-1}(x))}
=\sum_{i=1}^\kp\frac{\omega_i(x)}{t-M_i^{-1}(x)}
\eea 
where, the $M_i^{-1}(x)$ are a complete set of 
inverse branches of $M_\kp(x)$ and $l_{\kp}$ is the 
leading coefficient of $M_{\kp}(t)$. Rearranging 
the above equation and taking the limit as 
$t\to M_i^{-1}(x)$ we find,
\bea
\omega_i(x) = \frac{W_{\kp-1}(M_i^{-1}(x))}
{M_\kp^{\pri}(M_i^{-1}(x))}.\nonumber
\eea
By comparing coefficients of $t^{\kp-1}$ in (1.1) it 
is clear that,
\bea
\sum_{j=1}^{\kp}\omega_j(x) = c\nonumber
\eea
where $c$ is the quotient of the leading coefficients 
of $W_{\kp-1}(x)$ and $M_{\kp}(x)$. So, having 
chosen a particular polynomial $W_{\kp-1}(x)$ we are in 
a position to define the transformed weight 
function in the following way. Let $E_0$ be the 
interval $[a,b]$, and $P_n(x)$ the polynomials 
orthogonal with respect to $w_0(x)$ over $E_0$ : 
\bea
h_n\delta_{nm} &=& \int_a^b P_n(x)P_m(x)w_0(x)dx 
\nonumber\\
&=& \frac{1}{c}\int_a^b P_n(x)P_m(x)
\left(\sum_{j=1}^\kp\omega_j(x)\right)w_0(x)dx
\nonumber\\
&=& \frac{1}{c}\sum_{j=1}^\kp\int_{A_j}
P_n(M_\kp(t))P_m(M_\kp(t))\frac{W_{\kp-1}(t)}
{M_\kp^{\pri}(t)}M_\kp^{\pri}(t)w_0(M_\kp(t))dt 
\nonumber\\
&=&
\frac{1}{c}\int_EP_n(M_\kp(t))P_m(M_\kp(t))
W_{\kp-1}(t)w_0(t)dt,\nonumber
\eea
where
\bea
E := \cup_{j=1}^\kp A_j,\nonumber
\eea
and,
\bea
M_j^{-1} : [a,b] \to A_j.\nonumber
\eea
So the new sequence of polynomials is orthogonal 
on the disjoint set $E$ with respect to the 
weight function
\bea
\wh{w}(t) = W_{\kp-1}(t)w_0(M_\kp(t)).\nonumber
\eea
It was shown in \cite{Ger:Van} that 
\bea
\wh{P}_{n\kp}(x) = P_n(M_\kp(x)) \qquad n=0,1,2,... 
\nonumber
\eea
This generates every $\kp^{{\rm th}}$ polynomial 
in the sequence. The polynomials $\wh{P}_{n\kp+j}(x) \;\; 1
\leq j 
\leq \kp-1$, are referred to as the intermediate 
polynomials. Formulas relating the intermediate 
polynomials 
and recurrence coefficients of the new sequence 
to known quantities of the old sequence are given in 
\cite{Ger:Van}.


Some explicit expressions for the Chebyshev 
polynomials on several intervals, along with many other important results
can be found in \cite{Chen-Lawrence}. This paper is 
organized as follows. In this section we recall the theta
function  representation of these 
polynomials which was derived in \cite{Chen-Lawrence}. We
also present some  facts concerning the polynomials that
will be referred to throughout this paper.

Section 2 is a study of two auxiliary polynomials 
$S_g(x;n)$ and $G_{g+1}(x;n)$, of degree $g$ 
and $g+1$ respectively. The auxiliary polynomials are 
fundamental to the theory of the generalization of  
Chebyshev polynomials under consideration. We also derive for the first 
time non-linear difference equations satisfied by 
the recurrence coefficients that appear in the three 
term recurrence relation.

In section 3 we derive a product representation for 
our  polynomials from which 
one can obtain the coefficients of the polynomial in 
terms of the recurrence coefficients and the branch 
points. 
This representation is a new one that involves the 
auxiliary polynomials studied in section 2. 
It is particularly useful when we consider certain 
configurations of the branch points where 
the recurrence coefficients simplify into an algebraic form.

In section $4$ of this paper we show how certain polynomial 
mappings give rise to a countable number 
of special cases of our generalization of the  Chebyshev polynomials. 
These cases correspond to when the recurrence
coefficients are periodic. We give explicit 
formulas for the polynomials of degree not equal to an 
integer multiple of the period.  

The monic form of the Chebyshev polynomials on several intervals 
will be denoted by $P_n(x)$.  They  
are orthogonal on
$E:=[-1,\al]\cup_{j=1}^{g-1}[\bt_j,\al_{j+1}]
\cup[\bt_g,1]$ with respect to the following weight,
\bea
w(x) =  \frac{1}{\pi}
\sqrt{\frac{\prod_{j=1}^g(x-\al_j)}{(1-x^2)
\prod_{j=1}^g(x-\bt_j)}},
\eea
so that,
\bea
\int_EP_n(x)P_m(x)w(x)dx = \delta_{nm}h_n.
\nonumber
\eea
The polynomials of the second kind are 
\bea
Q_n(x) = \int_E\frac{P_n(x)-P_n(t)}{x-t}w(t)dt.
\eea
Both $P_n(x)$ and $Q_n(x)$ satisfy the recurrence relation
\bea
xu_n(x) = u_{n+1}(x) + b_{n+1}u_n(x) + a_nu_{n-1}(x),
\eea
see \cite{Akh1}. The Stieltjes transform of the weight function  is
\bea
\psi(z) := \int_E \frac{w(t)dt}{z-t} = 
\sqrt{\frac{\prod_{j=1}^g(z-\al_j)}{(z^2-1)
\prod_{j=1}^g(z-\bt_j)}},
\quad z\notin E.\nonumber
\eea
Note that, 
\bea
\psi(x) = \frac{\pi}{i}w(x).\nonumber
\eea
In \cite{Chen-Lawrence}, Chen and 
Lawrence considered a hyperelliptic Riemann surface of genus
$g$, 
\bea
y^2(z) = (z^2-1)\prod_{j=1}^g(z-\al_j)(z-\bt_j),\nonumber
\eea
from which a representation for the polynomials in 
terms of the associated Riemann theta function was 
constructed. They also gave a representation for 
the recurrence coefficients in terms of the theta 
functions. To formulate these representations the 
following facts regarding this Riemann surface are required. 

A canonical basis of cycles is chosen on the surface 
that is composed of $2g$ closed contours. We 
denote by the a-cycles the $g$ closed contours that 
lie on one sheet of the surface, and the b-cycles 
to be the remaining $g$ contours that travel over 
both sheets. For a surface of genus $g$ a basis of 
holomorphic differentials can be written in terms of 
$g$ linearly independent differentials,
\bea
\left\{\frac{dx}{y(x)},\frac{x\;dx}{y(x)},....,
\frac{x^{g-1}dx}{y(x)}\right\}.\nonumber
\eea
If we let,
\bea
A_{jk} = \int_{a_k}\frac{x^{g-j}}{y(x)}dx\nonumber
\eea
then the $d\omega_j$ defined as,
\bea
d\omega_j = \sum_{j=1}^g(A^{-1})_{jk}\frac{x^{g-j}}
{y(x)}dx,\nonumber
\eea
form an orthonormal basis such that,
\bea
\int_{a_j}d\omega_k = \delta_{jk}.\nonumber
\eea
We can then define the period matrix $B$ whose entries are,
\bea
B_{jk} = \int_{b_j}d\omega_k. \nonumber
\eea
The Riemann theta function associated with the surface 
is defined to be,
\bea
\vt(s;B) = \vt(s) := \sum_{t \in \BZ^g}exp(i\pi(t,Bt) + 2\pi
i(t,s)),\nonumber
\eea
where $s$ is a $g$-vector, and $(u,v)$ denotes the standard scalar 
product of two vectors. Convergence of the  series is
assured since $\Im B$ is a positive definite matrix. The
theta function has the following  properties.  For $t \in
\BZ^g$,
\bea
&&\vt(-s) = \vt(s)\nonumber\\
&&\vt(s+t) = \vt(s)\nonumber\\
&&\vt(s+Bt) = e^{-i\pi[(t,Bt)+2(t,s)]}\vt(s)\nonumber
\eea
Both the polynomials and the recurrence coefficients are 
given below. The expressions are taken from 
\cite{Chen-Lawrence}. However, in view of lemma $10.2$  in
\cite{Chen-Lawrence} and the  quasi-periodicity of the theta
functions we have simplified the expressions slightly. 
\bea
P_n(x) &=& \frac{1}{2}\left(\ce_n(\pg_x)+\cet_n(\pg_x)
\right),
\eea
where, for $n\geq 1$,
\bea
\ce_n(\pg_x) = 2e^{-n(\chi_0+\Omega(\pg_x))}
\frac{\vt(\up)\vt(\ux+\nbh)}{\vt(\up-\nbh)\vt(\ux)}
\nonumber
\eea
and,
\bea
\cet_n(\pg_x) = 2e^{n(\Omega(\pg_x)-\chi_0)}\frac{\vt(\up)\vt(\ux-\nbh)}{\vt(\up-\nbh)\vt(\ux)}.\nonumber
\eea
Using this notation,
\bea
Q_n(x) = \frac{\psi(x)}{2}
\left(\cet_n(\pg_x)-\ce_n(\pg_x)\right).
\eea
Here, $\pg_x$ is a point on the top sheet 
of the Riemann surface that corresponds to a point $x$ in 
the complex plane and, 
\bea
\Omega(\pg_x) &=& \int_{\pg_1}^{\pg_x}d\Omega\nonumber\\ 
&:=&\int_{\pg_1}^{\pg_x}
\frac{t^g+\sum_{j=0}^{g-1}k_jt^j}{y(t)}dt.\nonumber
\eea
The $k_j$ can be determined from the condition that,
\bea
\int_{a_j}d\Omega = 0, \;\;\; j=1, ...,  g.\nonumber
\eea
$\ux$ is defined to be a g-vector,
\bea
\ux := \left(\int_{p_1}^{p_x}d\omega_1, ...., 
\int_{p_1}^{p_x}d\omega_g\right)^T.\nonumber
\eea
and ${\bf u^{\pm}} := {\rm {\bf u}}_{\infty^{\pm}}.$ 
If $\hat{{\bf B}}$ is the vector whose components are,
\bea
\hat{{\bf B}}_j=\frac{1}{2\pi i}\int_{b_j}d\Omega, 
\;\;\; j=1,...,g. \nonumber
\eea
then,
\bea
\hat{{\bf B}} \equiv \up -\um \;\; \mbox{(mod $\BZ^g$)}.
\nonumber
\eea 
Finally, the equilibrium potential of the set $E$ is  
\bea
\chi_0 = \int_{\pg_1}^{\infty^{+}}
\left(\frac{t^g+\sum_{j=0}^gk_jt^j}{y(t)}
- \frac{1}{t}\right)dt,
\eea
With this notation, we have the following expressions 
for the recurrence coefficients,
\bea
a_n = \left\{ \begin{array}{ll}
2e^{-2\chi_0}\frac{\vt(\up+\hat{{\bf B}})}
{\vt(\up-\hat{{\bf B}})} & \textrm{if $n=1$}\\
e^{-2\chi_0}\frac{\vt(\up-(n-1)\hat{{\bf B}})
\vt(\up+n\hat{{\bf B}})}{\vt(\up+(n-1)\hat{{\bf
B}})\vt(\up-n\hat{{\bf B}})} & \textrm{if $n>1$}
\end{array} \right.
\eea
and,
\bea
b_n = \frac{1}{2}\sum_{j=1}^g(\bt_j-\al_j)
+
\sum_{j=1}^g(A^{-1})_{j1}
\left[{\scriptstyle 
\frac{\vt_j^{\pri}(\up-\hat{{\bf B}})}{\vt(\up-\hat{{\bf
B}})}-\frac{\vt_j^{\pri}(\up)}{\vt(\up)} +
\frac{\vt_j^{\pri}(\up-(n-1)\hat{{\bf
B}})}{\vt(\up-(n-1)\hat{{\bf
B}})}-\frac{\vt_j^{\pri}(\up-\nbh)}
{\vt(\up-\nbh)}}\right].
\eea
Here,
\bea
\vt_j^{\pri}({\bf u}) :=\frac{\partial}
{\partial u_j}\vt({\bf u})\nonumber
\eea 
Also,
\bea
h_n = 2e^{-2n\chi_0}\frac{\vt(\up+n{\bf \hat{B}})}
{\vt(\up-n{\bf \hat{B}})}.
\eea
\begin{theorem}
The polynomials of the second kind $Q_n(x)$ 
satisfy the following orthogonality relation,
\bea
\int_EQ_n(x)Q_m(x)\frac{1}{w(x)}dx = 
\delta_{nm}\tilde{h}_n,
\eea
where,
\bea
\tilde{h}_n = \pi^2h_n.\nonumber
\eea
\end{theorem}
{\it Proof :} 
We wish to evaluate the following integral,
\bea
I_j = \int_E\frac{Q_n(t)}{w(t)}t^jdt \qquad 0 
\leq j < n-1.\nonumber
\eea
We have,
\bea
Q_n(z) &=& P_n(z)\psi(z)
- \int_E\frac{P_n(t)}{z-t}w(t)dt\nonumber\\ 
&=&
P_n(z)\psi(z)-\frac{1}{z}\int_E
\frac{P_n(t)}{1-\frac{t}{z}}w(t)dt\nonumber\\
&=&
P_n(z)\psi(z)-\frac{1}{z}\int_EP_n(t)
\sum_{j=0}^{\infty}\left(\frac{t}{z}\right)^jdt
\nonumber\\
&=&
P_n(z)\psi(z)-\left(\frac{h_n}{z^{n+1}}+
O \left(\frac{1}{z^{n+2}}\right)\right)\nonumber
\eea
where the last line follows from the 
orthogonality of the $P_n(x)$. Therefore,
\bea
\frac{Q_n(z)}{\psi(z)} = P_n(z)-\frac{1}
{\psi(z)}\left(\frac{h_n}{z^{n+1}}
+ O\left(\frac{1}{z^{n+2}}\right)\right).\nonumber
\eea
We rewrite $I_j$ as a contour integral 
using a closed contour $\Lambda$ \cite{Chen-Lawrence} 
[Fig.
1],  in the slit complex plane that encircles all 
of the
branch points. By continuous deformation of the  contour
onto the intervals that make up $E$ we have,
\bea
I_j &=& \int_E\frac{Q_n(t)}{w(t)}t^jdt \nonumber\\
&=& \frac{\pi}{2i}\int_{\Lambda}
\frac{Q_n(z)}{\psi(z)}z^jdz.\nonumber
\eea
So $I_j$ can be written as,
\bea
I_j = \int_{\Lambda}P_n(z)z^jdz - 
\int_{\Lambda}\frac{\frac{h_n}{z^{n+1}} +
O\left(\frac{1}{z^{n+2}}\right)}{\psi(z)}z^jdz.
\nonumber
\eea
The differential $P_n(z)z^jdz$ has no residue at 
infinity for any $j$. The differential
$\frac{z^j}{z^{n+1}\psi(z)}dz$ has no residue at infinity
for $j=0,...,n-2$. Consequently, for these values of $j$,
both integrals are zero which proves the orthogonality of
the $Q_n(z)$ with respect to $\frac{1}{w(z)}$. The integral
$I_{n-1}$ is equal to the square of the weighted $L_2$ norm
of $Q_n(z)$. In this case there is a residue at infinity.
The residue is $-h_n$ and therefore,
\bea
\tilde{h}_n = \pi^2 h_n\nonumber
\eea 

$\Box$

\begin{theorem}

Let $\tilde{E}$ be the reflection of $E$ around the 
point $x=0$, namely 
\bea
\tilde{E}:=[-1,-\bt_g]\cup_{j=g-1}^{1}[-\al_{j+1},
-\bt_j]\cup[-\al_1,1].\nonumber 
\eea
Let $\{\tilde{P}_n(x)\}$ be the monic polynomials 
orthogonal on $\tilde{E}$ with respect to $\tilde{w}(x)$,
where,
\bea
\tilde{w}(x) = \frac{1}{\pi}
\sqrt{\frac{\prod_{j=1}^g(x+\bt_j)}{(1-x^2)
\prod_{j=1}^g(x+\al_j)}}.\nonumber
\eea
 Then we have for $n\geq g$
\bea
\tilde{P}_n(x) = \frac{(-1)^{n+g}}
{D_n\prod_{j=1}^g(x+\bt_j)}\left|
\begin{array}{cccccccc}
P_{n-g}(\bt_1) & \dots & \dots & \dots & P_{n+g}(\bt_1)  
\\
\vdots & \vdots & \vdots & \vdots & \vdots   \\
P_{n-g}(\bt_g) & \dots & \dots & \dots & P_{n+g}(\bt_g) 
 \\
Q_{n-g}(\al_1) & \dots & \dots & \dots & Q_{n+g}(\al_1)  
\\
\vdots & \vdots & \vdots & \vdots & \vdots  \\
Q_{n-g}(\al_g) & \dots & \dots & \dots & Q_{n+g}(\al_g)  
\\ 
P_{n-g}(-x) & \dots & P_n(-x) & \dots & P_{n+g}(-x) 
\end{array}
\right|.
\eea

If $n < g$ then,

\bea
\tilde{P}_n(x) = \frac{(-1)^{n+g}}
{D_n\prod_{j=1}^g(x+\bt_j)}\left|
\begin{array}{cccccccc}
0 & \dots & \dots & 0 & P_0(\bt_1) & 
\dots & P_{n+g}(\bt_1) \\
\vdots & \vdots & \vdots & \vdots 
& \vdots & \vdots & \vdots \\
0 & \dots & \dots & 0 & P_0(\bt_g) & \dots & 
P_{n+g}(\bt_g) \\
\al_1^0 & \al_1 & \dots & \al_1^{g-n-1} & 
Q_0(\al_1) & \dots & Q_{n+g}(\al_1)\\
\vdots & \vdots & \vdots & \vdots & \vdots & \vdots & 
\vdots \\
\al_g^0 & \al_g & \dots & \al_g^{g-n-1} & Q_0(\al_g) & 
\dots & Q_{n+g}(\al_g)\\
0 & 0 & 0 & 0 & P_0(-x) & \dots & P_{n+g}(-x) 
\nonumber
\end{array}
\right|.
\eea

In both equations, $D_n$ is equal to the coefficient 
of $P_{n+g}(-x)$ appearing in each determinant. 

\end{theorem}
{\it Proof:}
Starting with the orthogonality relation of the 
$\tilde{P}_n(x)$ we make the change of variables $x\to -x$.

\bea
\tilde{h}_n\delta_{nm} &=& \int_{\tilde{E}}
\tilde{P}_n(x)\tilde{P}_m(x)
\sqrt{\frac{\prod_{j=1}^g(x+\bt_j)}{(1-x^2)\prod_{j=1}^g
(x+\al_j)}}dx\nonumber\\
&=& \int_E\tilde{P}_n(-x)
\tilde{P}_m(-x)
\sqrt{\frac{\prod_{j=1}^g(x-\bt_j)}{(1-x^2)
\prod_{j=1}^g(x-\al_j)}}dx\nonumber\\
&=&
\int_E\tilde{P}_n(-x)\tilde{P}_m(-x)
\sqrt{\frac{\prod_{j=1}^g(x-\al_j)}
{(1-x^2)\prod_{j=1}^g(x-\bt_j)}} 
\left(\frac{\prod_{j=1}^g
(x-\bt_j)}{\prod_{j=1}^g(x-\al_j)}\right) dx
\nonumber
\eea

So the $\tilde{P}_n(-x)$ are orthogonal 
with respect to a weight function that is 
$w(x)$ multiplied
by a rational function. In general if we know the
polynomials orthogonal with respect to a weight $\mu(x)$ and
we wish to find the polynomials orthogonal with respect to
$R(x)\mu(x)$ where $R(x)$ is a rational function that is
positive on the interval of orthogonality we may invoke a
theorem of Uvarov that gives these polynomials in terms of
the original ones. The theorem was first published in
\cite{Uva} and also appears in \cite{Uva1}. In our case the
rational function has zeros and poles at the end points of
the subintervals that make up $E$. However, a brief
examination of the proof will show that the theorem is still
applicable in this case. The result follows directly from
the theorem. 
$\Box$

In section 2 we wish to study the following polynomials,
\bea
S_g(x;n):=P_n^2(x)\prod_{j=1}^g(x-\al_j)
- Q_n^2(x)(x^2-1)\prod_{j=1}^g(x-\bt_j)
\eea
and,
\bea
G_{g+1}(x;n):=P_n(x)P_{n-1}(x)
\prod_{j=1}^g(x-\al_j)-Q_n(x)Q_{n-1}(x)(x^2-1)
\prod_{j=1}^g(x-\bt_j).
\eea
These functions are fundamental in the theory of 
the generalizations of  Chebyshev polynomials considered here. 

\begin{theorem}

The functions $S_g(x;n)$ and $G_{g+1}(x;n)$ 
are polynomials of degree $g$ and $g+1$ respectively.

\end{theorem}
{\it Proof:} From (1.3) we know that,
\bea
Q_n(z) = P_n(z)\psi(z)-\int_E
\frac{P_n(t)}{z-t}w(t)dt.
\nonumber
\eea
Expanding the integrand in powers of $t$ we have,
\bea
\frac{1}{z-t} = \frac{1}{z}+\frac{t}{z^2}
+ \frac{t^2}{z^3}+....\nonumber
\eea    
Therefore, as $z \to \infty^{+}$, we have
\bea
P_n(z)\psi(z)-Q_n(z) = O
\left(\frac{1}{z^{n+1}}\right),
\nonumber
\eea
and similarly,
\bea
P_n(z)\psi(z)+Q_n(z) = O\left(z^{n-1}\right).
\nonumber
\eea
Multiplying the two together gives 
\bea
P_n^2(z)\psi^2(z)-Q_n^2(z)=O\left(\frac{1}{z^2}\right).
\nonumber
\eea
Therefore,
\bea
\prod_{j=1}^g(z-\al_j)P_n^2(z)-(z^2-1)
\prod_{j=1}^g(z-\bt_j)Q_n^2(z) = O\left(z^g\right).
\nonumber
\eea
Since the left hand side is a polynomial it follows that 
it must be of degree $g$. For (1.14) we form the following
product,
\bea
&&\left(P_n(z)\psi(z)-Q_n(z)\right)
\left(P_{n-1}(z)\psi(z)+Q_{n-1}(z)\right) \nonumber\\ 
&& \qquad \qquad \qquad = 
P_n(z)P_{n-1}(z)\psi(z)-Q_n(z)Q_{n-1}(z)-\psi(z)h_{n-1}
\nonumber
\eea
where we have used,
\bea
P_{n-1}(z)Q_n(z)-P_n(z)Q_{n-1}(z) = h_{n-1},
\nonumber
\eea
which is the Wronskian. As a result,
\bea
\prod_{j=1}^g(z-\al_j)P_n(z)P_{n-1}(z)-(z^2-1)
\prod_{j=1}^g(z-\bt_j)Q_n(z)Q_{n-1}(z) =
O\left(z^{g+1}\right).\nonumber
\eea
$\Box$

From (1.13) and (1.14) we have,
\bea
\ce_n(\pg_x)\cet_n(\pg_x) = \frac{S_g(x;n)}
{\prod_{j=1}^g(x-\al_j)}.
\eea
Also,
\bea
\ce_n(\pg_x)\cet_{n-1}(\pg_x) &=& 
\frac{G_{g+1}(x;n)+h_{n-1}y(x)}{\prod_{j=1}^g(x-\al_j)}\\
\cet_n(\pg_x)\ce_{n-1}(\pg_x) &=& 
\frac{G_{g+1}(x;n)-h_{n-1}y(x)}{\prod_{j=1}^g(x-\al_j)}.
\eea

An alternative representation for the 
 polynomials of the first and second
kind was also given in \cite{Chen-Lawrence} resulting from a
coupled system of differential equations. 
\bea
P_n(z) &:=& \rho_n(z)\cos\Psi_n(z)
\eea
and,
\bea
Q_n(z) &:=& i\psi(z)\rho_n(z)\sin\Psi_n(z),
\eea
where,
\bea
\rho_n(z) = \sqrt{\frac{S_g(z;n)}
{\prod_{j=1}^g(z-\al_j)}}\nonumber
\eea
and,
\bea
\Psi_n(z) = i\int_1^z
\frac{(\sum_{j=0}^{g-1}c_j(n)t^j-nt^g
+  \frac{1}{2}\sum_{j=1}^g\frac{y_j(n)}
{t-\gamma_j(n)})}{y}dt. 
\nonumber
\eea
where the branch of the square root is 
chosen in such a way that $\rho_n(z) \to \sqrt{2h_n}$, as 
$\Re z \to \infty$ and $2h_n$ is the leading coefficient  of
$S_g(z;n)$, to be shown later.  Here the path of integration
is from $1$ to an arbitrary complex point $z$ with the
property 
$\Im z > 0$, and the path of integration is entirely  in the
upper half plane. Expressions for the  real values are then
obtained by the analytical continuation of those above
allowing 
$z\to x \in \mathbb{R}$.

In these expressions, the $\ga_j(n)$ are the zeroes 
of $S_g(x,n)$, and $y_i = y(\pg_x)|_{\pg_x=\ga_j(n)}$. 
The
$c_j(n)$ are the first $g+1$ expansion coefficients of
$\frac{d}{dx}\ln\ce_n(\pg_x)$ in terms of the following
basis,
\bea
\left\{\frac{1}{y},\frac{x}{y},...,\frac{x^g}{y},
\frac{1}{x-\al_1},...,\frac{1}{x-\al_g},
\frac{y+y_1}{y(x-\ga_1(n))},...,\frac{y+y_g}
{y(x-\gamma_g(n))}\right\}
\nonumber
\eea  
This is the basis for meromorphic functions on 
the Riemann surface defined by $y(x)$ whose zeroes must
include simple zeroes at $\pm \infty$ and whose poles must
be from amongst the set
${\pg_1,\pg_{-1},\pg_{\al_j},\pg_{\bt_j},\pg_{\ga_j}
:j=1,..,g}$,
where the poles at the $\alpha$ points can be at most double
poles and the other poles must be simple poles. 

\setcounter{equation}{0}
\setcounter{theorem}{0}

\section{Evaluating the polynomials $S_g(x;n)$ and 
$G_{g+1}(x;n)$}

In this section we give an algorithm for determining 
the coefficients of $S_g(x;n)$ and $G_{g+1}(x;n)$ in terms
of the recurrence coefficients and the branch points. Let
\bea
S_g(x;n) = \sum_{j=0}^g\eta_j(n)x^j\nonumber
\eea
and,
\bea
G_{g+1}(x;n) = \sum_{j=0}^{g+1}\xi_j(n)x^j. 
\nonumber
\eea
Then we have,
\bea
\psi(z)P_n^2(z)-\frac{Q_n^2(z)}{\psi(z)} =
\frac{\eta_g(n)z^g+ ...+\eta_0(n)} {y(z)}
\eea
\bea
\psi(z)P_{n}(z)P_{n-1}(z) -  
\frac{Q_n(z)Q_{n-1}(z)}{\psi(z)}
=\frac{\xi_{g+1}(n)z^{g+1}+...+\xi_0(n)}{y(z)},
\eea
To find the $\eta_j(n)$ coefficients we 
integrate both sides of (2.1) along the contour $\Lambda$.
Note that,
\bea
\psi(z)=\frac{\pi}{i}W(z),\nonumber
\eea
where
\bea
W(z)=\frac{i}{\pi}
{\sqrt {\frac{\prod_{j=1}^{g}(z-\al_j)}{(z^2-1)
\prod_{j=1}^{g}(z-\bt_j)}}},
\nonumber
\eea
and for $t\in E,$
\bea
\lim_{\epsilon\to 0+}W(t+i\epsilon) =  w(t):=\frac{1}{\pi}
{\sqrt {\frac{\prod_{j=1}^{g}(t-\al_j)}
{(1-t^2)\prod_{j=1}^{g}
(t-\bt_j)}}}.
\eea
We obtain the following,
\bea
\frac{\pi}{i}\int_{\La}W(z)P_n^2(z)dz
- \frac{i}{\pi}\int_{\La}
\frac{Q_n^2(z)}{W(z)}dz =  \sum_{j=0}^{g}\eta_j(n)I_{j},
\eea
where,
\bea
I_j = \int_{\Lambda}\frac{z^j}{y(z)}dz. 
\nonumber
\eea
The $I_j$ can be evaluated by evaluating 
the residue at infinity. in which case,
\bea
I_j = 0 \qquad j < g \nonumber
\eea
and,
\bea
I_{g+k} = \int_{\Lambda}\frac{z^{g+k}}{\left[(z^2-1)
\prod_{j=1}^g(z-\al_j)(z-\bt_j)\right]^{1/2}}dz
\nonumber
\eea
which under the substitution $z=1/Z$ becomes,
\bea
I_{g+k} = -\int_{\Omega}\frac{dZ}{Z^{k+1}\left[(1-Z^2)\prod_{j=1}^g(1-\al_jZ)(1-\bt_jZ)\right]^{1/2}}.\nonumber
\eea
So in general the integrand of $I_{g+k}$ has a pole of order $k+1$ at infinity. Hence if,
\bea
f(Z) = \left[(1-Z^2)\prod_{j=1}^g(1-\al_jZ)(1-\bt_jZ)\right]^{-1/2}\nonumber
\eea
then,
\bea
I_{g+k} = -2\pi i\frac{f^{(k)}(0)}{k!}
\eea
where $f^{(k)}(Z)$ denotes the $k$th derivative of $f$. 
Deforming $\Lambda$ onto $E$ (2.4) becomes,
\bea
4h_n = \frac{i}{\pi}\eta_g(n)I_g.\nonumber
\eea
From (2.5) we have,
\bea
I_g = -2\pi i,\nonumber
\eea
and therefore,
\bea
\eta_g(n) = 2h_n.\nonumber
\eea
To find the remaining $\eta_j(n)$ coefficients 
we multiply (2.1) by successive powers of $z$ and perform
the same integration. This leaves us with the following
system of equations that we can use to solve for the
$\eta_g(n)$,
\bea
r_0(n) &=& I_g\eta_g(n)\nonumber\\ 
r_1(n) &=& I_{g+1}\eta_g(n) + I_g\eta_{g-1}(n)\nonumber\\
r_2(n) &=& I_{g+2}\eta_g(n)
+
I_{g+1}\eta_{g-1}(n)+I_g\eta_{g-2}(n)
\nonumber\\ 
&\vdots&
\nonumber\\ 
r_g(n) &=& \sum_{j=0}^gI_{g+j}\eta_{j}.
\nonumber
\eea
Here,
\bea
r_k(n) = \frac{4\pi}{i}
\int_Ez^kP_n^2(z)w(z)dz .
\eea
These constants are obtained by iterating the 
three term recurrence relation $k$ times in (2.6). The
recurrence relation can also be expressed in matrix form by,
\bea
Lp = zp,\nonumber
\eea
where, 
\bea
L = \left(
\begin{array}{cccccccc}
b_1 & 1 & 0 & 0 & 0 & 0 & 0  \\
a_1 & b_2 & 1 & 0 & 0 & 0 & 0  \\
0 & a_2 & b_3 & 1 & 0 & 0 & 0  \\
0 & 0 & a_3 & b_4 & 1 & 0 & 0  \\
0 & 0 & 0 & a_4 & b_5 & 1 & 0  \\
0 & 0 & 0 & 0 & \ddots & \ddots & \ddots  \\ 
0 & 0 & 0 & 0 & 0 &\ddots & \ddots  \nonumber
\end{array}
\right)
\;\;\;\;\;\;\;\;\;\;\;\;
p = \left(
\begin{array}{cccccccc}
P_0(z) \\
P_1(z) \\
P_2(z) \\
P_3(z) \\
\vdots \\
\vdots \\
\vdots \nonumber
\end{array}
\right).
\eea  
In this case, 
\bea
z^kp &=& L^kp\nonumber\\
zP_n(z) &=& \sum_{m=0}^{\infty}L_{nm}p_m \nonumber\\
z^kP_n(z) &=& \sum_{m=0}^{\infty}[L^k]_{nm}P_m(z).
\nonumber
\eea
From the expression for $r_k(n)$ in (2.6) it is 
clear that,
\bea
r_k(n) = \frac{4\pi h_n}{i}[L^k]_{n-1,n}. 
\eea
To find the $\xi_j(n)$ we use the same method. 
We know from the proof of theorem (1.3) that,
\bea
\xi_{g+1}(n) = h_{n-1}.\nonumber
\eea 
Multiplying (2.2) by successive powers of $z$ and 
integrating we have the following system of equations,
\bea
0 &=& I_{g+1}\xi_{g+1}(n)+I_g\xi_g(n)
\nonumber\\
s_1(n) &=& I_{g+2}\xi_{g+1}(n)+I_{g+1}\xi_g(n) +
I_g\xi_{g-1}(n)\nonumber\\ s_2(n) &=&
I_{g+3}\xi_{g+1}(n)+I_{g+2}\xi_g(n)
+ I_{g+1}\xi_{g-1}(n)+I_g\xi_{g-2}(n)
\nonumber\\
&\vdots& \nonumber\\ s_g(n) &=&
\sum_{j=0}^{g+1}I_{g+j}\xi_j(n).
\nonumber
\eea
Here,
\bea
s_k(n) = \frac{4\pi}{i}\int_Ez^kP_n(z)P_{n-1}(z)w(z)dz.
\nonumber
\eea
In terms of the Jacobi matrix,
\bea
s_k(n) = \frac{4\pi h_n}{i}[L^k]_{nn}.
\nonumber
\eea
For $g=1$ we have,
\bea
S_1(x,n) &=& 2h_n\left(x+b_{n+1}-\frac{(\al+\bt)}{2}\right)
 \qquad n\geq 1\nonumber\\
G_2(x,n) &=& h_{n-1}
\left(x^2-\frac{\al+\bt}{2}x  +
2a_n-\frac{1}{8}(\al-\bt)^2 -  \frac{1}{2}\right).
\qquad n\geq 2\nonumber
\eea
\bea
S_1(x,0) &=& (x-\al)\nonumber\\
G_2(x,1) &=& (x-b_1)(x-\al).\nonumber
\eea
The auxiliary polynomials for $g=2$ can be found 
in appendix B.

{\bf Non-linear difference equations}

We can use the polynomials studied in this section 
to derive a pair of non linear difference equations for each
genus that are satisfied by the recurrence coefficients
$\{a_n\}$ and $\{b_n\}$. Using the three term recurrence
relation, 
\bea
\frac{S_g(z;n)}{\prod_{j=1}^g(z-\al_j)} 
&=& \ce_n(\pg_z)\cet_n(\pg_z)\nonumber\\
&=& (z-b_n)^2\ce_{n-1}(\pg_z)\cet_{n-1}(\pg_z)
\nonumber\\
&-&a_{n-1}(z-b_n)(\ce_{n-2}(\pg_z)\cet_{n-1}
(\pg_z)+\cet_{n-2}(\pg_z)\ce_{n-1}(\pg_z))
+ a_{n-1}^2\ce_{n-2}(\pg_z)\cet_{n-2} (\pg_z)
\nonumber\\
&=&\frac{(z-b_n)^2S_g(z;n-1)+a_{n-2}^2
S_g(z;n-2)-2a_{n-1}(z-b_n)G_{g+1}(z;n-1)}
{\prod_{j=1}^g(z-\al_j)}.
\eea
Evaluating both sides at $z=b_n$ gives,
\bea
S_g(b_n;n) = a_{n-1}^2S_g(b_n,n-2).\nonumber 
\eea 
If we write this out fully for $g=1$ we have,
\bea
a_n\left(b_{n+1}+b_n-\frac{\al+\bt}{2}\right) 
= a_{n-1}\left(b_n+b_{n-1}-\frac{\al+\bt}{2}\right).
\eea
This is our first difference relation for the 
recurrence coefficients. To get the second we note that,
\bea
\frac{G_{g+1}(z;n)}{\prod_{j=1}^g(z-\al_j)} 
&=&
\frac{1}{2}\left(\ce_n(\pg_z)\cet_{n-1}(\pg_z)
+ \cet_n(\pg_z)\ce_{n-1}(\pg_z)\right)\nonumber\\
&=&
\frac{1}{2}\Big(\big((z-b_n)\ce_{n-1}
(\pg_z)-a_{n-1}\ce_{n-2}(\pg_z)\big)
\cet_{n-1}(\pg_z)\nonumber\\
&&\;\;\;\;+\big(z-b_n)\cet_{n-1}(\pg_z)-a_{n-1}
\cet_{n-2}(\pg_z)\big)\ce_{n-1}(\pg_z)\Big)\nonumber\\
&=&(z-b_n)\ce_{n-1}(\pg_z)\cet_{n-1}(\pg_z)
- \frac{a_{n-1}}{2}
\left(\ce_{n-2}(\pg_z)\cet_{n-1}(\pg_z)
+ \cet_{n-2}(\pg_z)\ce_{n-1}(\pg_z)\right)
\nonumber\\
&=&
\frac{(z-b_n)S_g(z;n-1)-a_{n-1}G_{g+1}(z;n-1)}
{\prod_{j=1}^g(z-\al_j)}. 
\eea
If we evaluate this at $z=b_n$ we get ,
\bea
G_{g+1}(b_n;n) = -a_{n-1}G_{g+1}(b_n;n-1).
\nonumber
\eea
Writing this out for genus 1,
\bea
a_n+a_{n-1} = \frac{1}{2}+\frac{(\bt-\al)^2}{8}
+ \frac{\al+\bt}{2}b_n-b_n^2.
\eea
This is our second difference equation for the 
recurrence coefficients. Both equations are valid for $n
\geq 3$.

\setcounter{equation}{0}
\setcounter{theorem}{0}

\section{Product representation for $P_n(x)$ and $Q_n(x)$} 

In this section we use the polynomials derived in section $2$ 
to derive a representation for the generalized Chebyshev polynomials 
where the coefficients of the powers of $x$ are given in terms of 
the branch points and the recurrence coefficients. From (1.5) we have
\bea
P_n(x) = \frac{1}{2}\left(\ce_n(\pg_x)+\cet_n(\pg_x)\right).\nonumber
\eea
We define the following functions,
\bea
f_+(x;n) := \frac{\ce_n(\pg_x)}{\ce_{n-1}(\pg_x)}\nonumber\\
f_-(x;n) := \frac{\cet_n(\pg_x)}{\cet_{n-1}(\pg_x)}.\nonumber
\eea
Using (1.16) and (1.17) we can write,
\bea
f_+(x;n) &=& \frac{\ce_n(\pg_x)\cet_{n-1}(\pg_x)}
{\ce_{n-1}(\pg_x)\cet_{n-1}(\pg_x)} = \frac{G_{g+1}(x;n)+h_{n-1}y(x)}{S_g(x;n-1)}. 
\eea
Similarly,
\bea
f_-(x;n) &=& \frac{\cet_n(\pg_x)\ce_{n-1}(\pg_x)}
{\cet_{n-1}(\pg_x)\ce_{n-1}(\pg_x)} = \frac{G_{g+1}(x;n)-h_{n-1}y(x)}{S_g(x;n-1)}.
\eea
Since, $\ce_0(\pg_x) = 1$,
we have, for any $x\in \mathbb{C}$  
\bea
P_n(x) = \frac{1}{2}\left(\prod_{j=1}^nf_+(x,j)+\prod_{j=1}^nf_-(x,j)\right).
\eea
Similarly,
\bea
Q_n(x) = \frac{\psi(x)}{2}\left(\prod_{j=1}^nf_+(x,j)-\prod_{j=1}^nf_-(x,j)\right).
\eea
(3.3) and (3.4) are valid for any genus $g$ and require only 
knowledge of the branch points and the recurrence coefficients. 
They do not require any knowledge of the $\gamma$ points. If we look at the genus 
1 case we find that,
\bea
f_{\pm}(x;n) = \frac{x^2-\frac{\al+\bt}{2}x+2a_n 
- \frac{1}{8}(\al-\bt)^2-\frac{1}{2} \pm \sqrt{(x^2-1)(x-\al)(x-\bt)}}
{2(x+b_{n+1}-\frac{\al+\bt}{2})}.\nonumber
\eea  
If we take the limit as $\al \to \bt$ we see that,
\bea
S_1(x,n) \to 2h_n(x-\al) \qquad n \geq 1\nonumber
\eea
\bea
G_2(x,n) \to h_{n-1}x(x-\al),\nonumber
\eea
and,
\bea
f_{\pm}(x,n) \to \frac{1}{2}(x\pm \sqrt{x^2-1}) 
= \frac{1}{2}e^{\pm i\theta} \qquad x=\cos\theta \;\;\;\; n \geq 2.\nonumber
\eea
Note that,
\bea
f_{\pm}(x,1) \to e^{\pm i\theta}. \nonumber
\eea
Thus we can see that (3.3) is a natural generalization of the following formula 
for the classical Chebyshev 
polynomials,
\bea
T_n(x) = \frac{(x+\sqrt{x^2-1})^n+(x-\sqrt{x^2-1})^n}{2^n}.\nonumber
\eea
In general, (3.3) and (3.4) are very effective representations of the polynomials 
for the purpose of extracting the polynomial coefficients.
The appearance of the theta function of the Riemann surface is still apparent in 
the expression for the recurrence
coefficients. One may also wish to study certain special cases of the Akhiezer 
polynomials when the recurrence coefficients become algebraic functions of   
the branch points $\al_i$ and $\bt_i$. In this case (3.1) is particularly useful 
as one could insert directly the value of the recurrence coefficients. 
In the next section we will encounter some of these special cases.
The following formulas follow naturally from (3.1):
\bea
P_n(x) = \frac{G_{g+1}(x,n)}{S_g(x,n-1)}P_{n-1}(x)+\frac{h_{n-1}(x^2-1)\prod_{j=1}^g(x-\bt_j)}{S_g(x,n-1)}Q_{n-1}(x),
\eea
and the corresponding formula for the $Q_n$ is,
\bea
Q_n(x) = \frac{G_{g+1}(x,n)}{S_g(x,n-1)}Q_{n-1}(x)+\frac{h_{n-1}\prod_{j=1}^g(x-\al_j)}{S_g(x,n-1)}P_{n-1}(x).
\eea
Letting $\al_i \to \bt_i$ in (3.3) we obtain as expected the analogous formula for the Chebyshev case:
\bea
2T_n(x) = xT_{n-1}(x)+(x^2-1)U_{n-1}(x).\nonumber
\eea
Here $T_n(x)$ and $U_n(x)$ are monic Chebyshev polynomials of the first and second kind respectively.

Now we show how the product representations of (3.1) and (3.2) can be used to re-express the differential relations derived in \cite{Chen-Lawrence}, where $P_n(x)$ and the $Q_n(x)$ were shown to satisfy,
\bea
P_n^{\pri}(x) = f_1(x;n)P_n(x)+f_2(x;n)Q_n(x),
\eea
and
\bea
Q_n^{\pri}(x) = f_3(x;n)P_n(x)+f_4(x;n)Q_n(x),
\eea
where,
\bea
f_1(x;n) &=& \frac{1}{2}\sum_{j=1}^g\left(\frac{1}{x-\gamma_j}-\frac{1}{x-\al_j}\right)\nonumber
\eea
\bea
f_2(x;n) &=& \frac{nx^g-\sum_{j=0}^{g-1}c_jx^j-\frac{1}{2}\sum_{j=1}^g\frac{y_j}{x-\gamma_j}}{\prod_{j=1}^g(x-\al_j)},\nonumber
\eea
and,
\bea
f_3(x;n) &=& \psi^2(x)f_2(x)\nonumber
\eea
\bea
f_4(x;n) &=& \frac{\psi^{\pri}(x)}{\psi(x)}+f_1(x).\nonumber
\eea
The definitions of all the constants appearing in the above formulas can be found in section 1. From \cite{Chen-Lawrence},
\bea
f_1(x;n) = \frac{\rho_n^{\pri}(x)}{\rho_n(x)},\nonumber
\eea
and,
\bea
f_2(x;n) = i\frac{\Psi_n^{\pri}(x)}{\psi(x)}.
\eea
Using (3.1) we can express $f_1(x;n)$ and $f_2(x;n)$ in more explicit terms. Recall that,
\bea
\rho_n(x) = \sqrt{\frac{S_g(x,n)}{\prod_{j=1}^g(x-\al_j)}}.\nonumber
\eea
Therefore, we can write,
\bea
f_1(x;n) = \frac{\rho_n^{\pri}(x)}{\rho_n(x)}=\frac{1}{2}\left(\frac{S_g^{\pri}(x,n)}{S_g(x,n)}-\sum_{j=1}^g\frac{1}{x-\al_j}\right).
\eea
This is simply a restatement of the earlier definition if we recall that the $\gamma_i(n)$ are the zeroes of $S_g(x,n)$.  
In order to re-express $f_2(x;n)$, we observe that 
\bea
\prod_{j=1}^nf_{\pm}(x,j) = \rho_n(x)e^{\pm i\Psi_n(x)}.\nonumber
\eea
Therefore,
\bea
\prod_{j=1}^n\frac{f_+(x,j)}{f_-(x,j)} = e^{2i\Psi_n(x)} \quad {\rm or} 
\quad
\ln\prod_{j=1}^n\frac{f_+(x,j)}{f_-(x,j)} = 2i\Psi_n(x) \qquad ({\textrm {mod}}\; 2\pi).\nonumber
\eea
From the definitions of $f_{\pm}(x,j)$ we have,
\bea
2i\Psi_n(x) = \sum_{j=1}^n\ln\left(\frac{G_{g+1}(x;j)+h_{j-1}y(x)}{G_{g+1}(x;j)-h_{j-1}y(x)}\right) \;\; ({\textrm {mod}}\; 2\pi).
\eea
Hence,
\bea
i\Psi_n^{\pri}(x) = \sum_{j=1}^n\left[\frac{y^{\pri}(x)G_{g+1}(x;j)-y(x)G_{g+1}^{\pri}(x;j)}{S_g(x;j)S_g(x;j-1)}\right]h_{j-1}.
\eea

Substituting the above into (3.12) gives us the differential relations (3.7) and (3.8) completely in terms of the recurrence coefficients and the branch points. 

Note that, generally when $\ga_l(k) \in (\al_l,\bt_l) \;\; k \geq 1$,  we can evaluate $\Psi_n(\al_l)$ by taking the limit $x \to \al_l$ in (3.11). From (1.13) and (1.14) we see that for $x=\al_l$ the only non-zero term appearing in the sum in (3.12) is the $j=1$ term. For $j>1$, $G_{g+1}(\al_l;j)$ cannot be zero since this would imply, from (1.13) and (1.14) that $\ga_l(j) = \al_l$. We have,
\bea
\lim_{x \to \al_l} 2i\Psi_n(x) = \lim_{x \to \al_l} \ln \left(\frac{G_{g+1}(x;1)+y(x)}{G_{g+1}(x;1)-y(x)}\right) \qquad ({\textrm {mod}}\; 2\pi).\nonumber
\eea
We know that,
\bea
G_{g+1}(x;1) = (x-b_1)\prod_{j=1}^g(x-\al_j)\nonumber
\eea
and therefore,
\bea
\lim_{x \to \al_l} \Psi_n(x) = \frac{\pi}{2} \qquad ({\textrm {mod}}\; \pi).\nonumber
\eea

{\bf 3b. Discriminants}

In this section we derive an expression for the discriminant when $g=1$. Once we have knowledge of the differential relations satisfied by a sequence of polynomials one can in general use this information to say something about the discriminant. The discriminant is useful when we consider certain electrostatic problems regarding the zeroes of the polynomials.  See \cite{Szego} and 
\cite{Ism3} for more details. The expressions
 we derive in this section are restricted to the cases when the polynomials do not have zeroes at the branch points. If we write,
\bea
P_n(x) = \prod_{j=1}^n(x-x_{j,n})\nonumber
\eea
then the discriminant is defined to be the following,
\bea
D[P_n(x)] := \prod_{1\leq i < j\leq n}(x_{i,n}-x_{j,n})^2.
\eea
Since
\bea
P_n^{\pri}(x) = \sum_{j=1}^n\frac{P_n(x)}{(x-x_{j,n})}\nonumber
\eea
it is easily verified that,
\bea
(-1)^{\frac{n(n-1)}{2}}\prod_{j=1}^nP_n^{\pri}(x_{j,n}) = D[P_n(x)].\nonumber
\eea
From (3.7), assuming that $P_n(x)$ has no zeroes at the branch points,  we have,
\bea
P_n^{\pri}(x_{j,n}) = f_2(x_{j,n})Q_n(x_{j,n}).\nonumber
\eea
To evaluate $Q_n(x_{j,n})$, recall that,
\bea
P_n(x) = \rho_n(x)\cos\Psi_n(x),\nonumber
\eea
and
\bea
Q_n(x) = i\psi(x)\rho_n(x)\sin\Psi_n(x).\nonumber
\eea
Since we have,
\bea
\frac{1}{i\psi(x)}\frac{Q_n(x)}{P_n(x)} = \tan \Psi_n(x),\nonumber
\eea
it must be that,
\bea
\tan \Psi_n(x_{j,n}) = \pm \infty \nonumber,
\eea
and therefore if $x_{j,n}$ is a zero of $P_n(x)$ then we do have,
\bea
\cos \Psi_n(x_{j,n}) = 0. \nonumber
\eea
Hence,
\bea
\sin \Psi_n(x_{j,n}) = \pm 1 .\nonumber
\eea
Therefore, it is always true that,
\bea
|Q_n(x_{j,n})| = |i\psi(x_{j,n})\rho_n(x_{j,n})|.\nonumber
\eea
The sign depends on the value of $n$ and $j$. To show this dependence, we note that the following quantity,
\bea
i\psi(x_{j,n})\rho_n(x_{j,n}) = i\sqrt{\frac{2h_n}{x_{j,n}^2-1}\prod_{l=1}^g
\frac{x_{j,n}-\ga_l(n)}{x_{j,n}-\bt_l}},\nonumber
\eea
is positive. Since theorem (7.1) in \cite{Chen-Lawrence} states that any zero 
of $P_n(x)$ on the interval $(\al_l,\bt_l)$ must lie in $(\al_l,\ga_l)$ and 
since the zeroes of $P_n(x)$ and $Q_n(x)$ interlace, we must have  
\bea
\sin\Psi_{2n}(x_{j,2n}) = (-1)^{j}\nonumber
\eea
and, 
\bea
\sin\Psi_{2n-1}(x_{j,2n-1}) = (-1)^{j+1}.\nonumber
\eea
Therefore,
\bea
Q_n(x_{j,n}) &=& (-1)^{n+j}i\psi(x_{j,n})\rho_n(x_{j,n}),\nonumber
\eea 
and, 
\bea
D[P_n(x)] &=& (-1)^{\frac{n(n-1)}{2}}\prod_{j=1}^nP_n^{\pri}(x_{j,n})\nonumber\\ 
&=& (-1)^{\frac{n(n-1)}{2}}\prod_{j=1}^n(-1)^{n+j+\frac{1}{2}}\sqrt{2h_n}\sqrt{\frac{\prod_{i=1}^g(x_{j,n}-\ga_i(n))}{(x_{j,n}^2-1)\prod_{i=1}^g(x_{j,n}-\bt_i)}}f_2(x_{j,n})\nonumber\\
&=& (-1)^{\frac{n}{2}}(2h_n)^{\frac{n}{2}}\sqrt{\frac{\prod_{j=1}^gP_n(\gamma_j(n))}{P_n(1)P_n(-1)\prod_{j=1}^gP_n(\bt_j)}}\prod_{j=1}^nf_2(x_{j,n}).
\eea

Now in this most general form we can not yet evaluate the discriminant for general $g$ since we do not know how to evaluate the product of the $f_2(x_{j,n})$. However, examining the form of $f_2(x;n)$ we see that it is a rational function. In certain special cases we may know enough about this function to factor it in which case we can express the unknown quantity in (3.17) as $P_n(x)$ evaluated at the zeroes and poles of $f_2(x;n)$.  For $g=1$ we can always do this, in which case,
\bea
D[P_n(x)]
=(-1)^{\frac{n}{2}}(2h_n)^{\frac{n}{2}}\sqrt{\frac{P_n(\gamma(n))}{P_n(1)P_n(-1)P_n(\bt)}}\prod_{j=1}^nf_2(x_{j,n}).\nonumber
\eea
Now,
\bea
f_2(x;n) = \frac{n(x-r_+(n))(x-r_-(n))}{(x-\gamma(n))(x-\al)}\nonumber
\eea
where,
\bea
r_{\pm}(n)=\frac{(c_0/n+\gamma)\pm\sqrt{(c_0/n-\gamma)^2+2y_1/n}}{2}.\nonumber
\eea
Therefore, we have for the square of the discriminant
\bea
D^2[P_n(x)]=(2n^2h_n)^n\frac{(-1)^n P_n^2(r_+)P_n^2(r_-)}{P_n^2(\al)P_n(1)P_n(-1)P_n(\bt)P_n(\ga(n))}
\eea
The product $P_n(\bt)P_n(\ga_n)$ is always positive since we know from \cite{Chen-Lawrence} that $P_n(x)$ cannot 
have a zero in the interval $(\ga(n),\bt)$.  
A similar calculation for the polynomials of the second kind gives,
\bea
D^2[Q_n(x)] = (2n^2h_n)^{n-1}\frac{Q_n^2(r_+(n))Q_n^2(r_-(n))}{Q_n^2(1)Q_n^2(-1)Q_n^2(\bt)Q_n(\gamma(n))Q_n(\al)}.
\eea
Again, in this case we require that $Q_n(x)$ have no zeroes at the branch points. 

These discriminants are valid for genus 1. Indeed we know explicitly the values of $\ga(n)$ and $r_{\pm}(n)$ for $g=1$. It is clear from the derivations in this section, that the equivalent expressions for general genus $g$ are,
\bea
D^2[P_n(x)] = (2n^2h_n)^n\frac{(-1)^n\prod_{j=1}^{2g}P_n^2(r_j(n))}{P_n(1)P_n(-1)\prod_{j=1}^g[P_n(\bt_j)P_n(\ga_j(n))P_n^2(\al_j)]},
\eea
and,
\bea
D^2[Q_n(x)] = (2n^2h_n)^{n-1}\frac{\prod_{j=1}^{2g}Q_n^2(r_j(n))}{Q_n^2(1)Q_n^2(-1)\prod_{j=1}^g[Q_n^2(\bt_j)Q_n(\gamma_j(n))Q_n(\al_j)]},
\eea
where the $\{r_j(n)\}$ are the $2g$ roots of the polynomial,
\bea
\left[\prod_{j=1}^g(x-\ga_j(n))\right]\left(nx^g-\sum_{j=0}^{g-1}c_j(n)x^j-\frac{1}{2}\sum_{j=1}^g\frac{y_j}{x-\ga_j(n)}\right)\nonumber.
\eea

\setcounter{equation}{0}
\setcounter{theorem}{0}

\section{Connection with Polynomial mappings}

The theta function expressions for the recurrence coefficients $a_n$ and $b_n$ of the generalized Chebyshev polynomials are given in (1.5) and (1.6). Consider for the moment two different values of the recurrence coefficients, for example, $b_n$ and $b_{n+\kp}$, where $n$ and $\kp$ are both integers. If $\kp$ is chosen such that $\kp{\bf \hat{B}} \in \mathbb{Z}^g$ then from the periodicity properties of the theta functions we have,
\bea
b_{n+\kp} = b_n. \nonumber
\eea 
From the expression for the $a_n$ coefficients, we will also have,
\bea
a_{n+\kp} = a_n  \qquad n \geq 2.\nonumber
\eea
When the entries of the vector ${\bf \hat{B}}$ are all rational numbers, we can always find an integer $\kp$ such that $\kp{\bf \hat{B}} \in \mathbb{Z}^g$. If we take the smallest integer that satisfies this requirement we will then have the smallest integer over which the recurrence coefficients themselves repeat, that is the recurrence coefficients are periodic with period $\kp$. If at least one of the entries is an irrational number we are unable to find an integer $\kp$ that satisfies the periodicity requirement. 

We will study the effect that this periodicity has on our previous constructions and use this to explain the connection between the generalized Chebyshev polynomials and non-trivial polynomial mappings \cite{Ger:Van}. Throughout this section $\kp$ is a given fixed integer. First consider $P_{n\kp}(x)$. From 
(1.5), and with the periodicity condition, we see that the $\theta$ factors 
cancel in equation (1.5), therefore 
\bea
P_{n\kp}(x) = \frac{\Delta_{n\kp}}{2}(e^{n\kp(\Omega(p_x))}+e^{-n\kp(\Omega(p_x))}),
\eea
hence,
\bea
P_\kp(x) = \frac{\Delta_\kp}{2}(e^{\kp(\Omega(p_x))}+e^{-\kp(\Omega(p_x))}).
\eea
Here,
\bea
\Delta_{n\kp} = 2e^{-n\kp\chi_0}.
\eea
Therefore,
\bea
P_{n\kp}(x) &=& 2^{n-1}\Delta_{n\kp}T_n\left(\frac{e^{\kp\Omega(x)}+e^{-\kp\Omega(x)}}{2}\right)\nonumber\\
&=&\Delta_\kp^nT_n\left(\frac{P_\kp(x)}{\Delta_\kp}\right),
\eea
where $T_n(x)$ is the monic Chebyshev polynomial. From the expression for $h_n$ in (1.10) we see that,
\bea
h_\kp = 2e^{-2\kp\chi_0}.\nonumber
\eea
It follows that,
\bea
\Delta_\kp=\sqrt{2h_\kp},\nonumber
\eea
and therefore,
\bea
P_{n\kp}(x)=(2h_\kp)^{\frac{n}{2}}T_n\left(\frac{P_\kp(x)}{\sqrt{2h_\kp}}\right).
\eea
The corresponding formula for the polynomials of the second kind is,
\bea
Q_{n\kp}(x) = (2h_\kp)^\frac{n-1}{2}Q_\kp(x)U_n\left(\frac{P_\kp(x)}{\sqrt{2h_\kp}}\right),
\eea
where $U_n(x)$ is the Chebyshev polynomial of the second kind (degree $n-1$). 
Recall that,
\bea
\ce_n(\pg_x)\cet_n(\pg_x) = \frac{S_g(x;n)}{\prod_{j=1}^g(x-\al_j)}.\nonumber
\eea
When $\kp{\bf \hat{B}} \in \BZ^g$, we have,
\bea
\ce_n(\pg_x)\cet_n(\pg_x) = 4e^{-2n\chi_0}.\nonumber
\eea 
Therefore, since the $\ga_j(n)$ are the zeroes of $S_g(x;n)$ we must have
$\ga_j(\kp)=\al_j$. It will be shown later in this section that the maximum value of $g$ permitted
for the polynomials is $\kp-1$. Therefore, if we put $g=\kp-1$ and $n=\kp$ in (1.13) we see 
immediately that,
\bea
Q_\kp(x) = \prod_{j=1}^{\kp-1}(x-\al_j).\nonumber
\eea
Hence,
\bea
Q_{n\kp}(x) = (2h_\kp)^\frac{n-1}{2}U_n\left(\frac{P_\kp(x)}{\sqrt{2h_\kp}}\right)\prod_{j=1}^{\kp-1}(x-\al_j).\nonumber
\eea
In the language of (3.3), (4.5) becomes,
\bea
P_{n\kp}(x) = \frac{(\prod_{j=1}^\kp f_+(x,j))^n+(\prod_{j=1}^\kp f_-(x,j))^n}{2^{n}},\nonumber
\eea
and (4.6) becomes,
\bea
Q_{n\kp}(x) = \psi(x)\frac{(\prod_{j=1}^\kp f_+(x,j))^n-(\prod_{j=1}^\kp f_-(x,j))^n}{2^n}.\nonumber
\eea
Using the theta function expressions, we can verify the following,
\bea
f_+(x,\kp+1) = \frac{\cet_{\kp+1}(\pg_x)}{\cet_\kp(\pg_x)} = \frac{\tilde{\Delta}_{\kp+1}}{\tilde{\Delta}_\kp\tilde{\Delta}_1}f_+(x,1).\nonumber
\eea
Similarly we can deduce that,
\bea
\frac{\tilde{\Delta}_{\kp+1}}{\tilde{\Delta}_\kp\tilde{\Delta}_1} = \frac{1}{2},\nonumber
\eea
and therefore,
\bea
f_{\pm}(x,\kp+1) = \frac{1}{2}f_{\pm}(x,1).\nonumber
\eea
An analogous calculation shows that,
\bea
f_{\pm}(x,\kp+j) = f_{\pm}(x,j), \;\;\;\; 1<j<\kp,\nonumber
\eea
as we expect from looking at the explicit structure of the $f_{\pm}(x,j)$.

If we know the period $\kp$ then we can use the $P_\kp(x)$ with (4.5) to generate every $\kp^{th}$ polynomial in the sequence. To obtain the intermediate polynomials $P_{n\kp+j}$ for $j=1..\kp-1$ we use the fact that
\bea
\ce_{n\kp+j}(\pg_x) = \frac{1}{2}\left(\prod_{l=1}^jf_+(x;l)\right)\ce_{n\kp}(\pg_x) = \frac{1}{2}\ce_j(\pg_x)\ce_{n\kp}(\pg_x),\nonumber
\eea
and therefore,
\bea
P_{n\kp+j}(x) = \frac{1}{2}\left(P_j(x)P_{n\kp}(x)+\frac{1}{\psi^2(x)}Q_j(x)Q_{n\kp}(x)\right), \qquad  1 \leq j < \kp , \;\; n \geq 1
\eea
Similarly,
\bea
Q_{n\kp+j}(x) = \frac{1}{2}\left(P_j(x)Q_{n\kp}(x)+Q_j(x)P_{n\kp}(x)\right), \qquad 1 \leq j < \kp, \;\; n \geq 1
\eea
So if we know the first $\kp$ polynomials of the first and second kind we can generate all of the remaining polynomials using (4.5), (4.6), (4.7) and (4.8).

From the theory of the classical Chebyshev polynomials we have,
\bea
T_m(x)T_n(x) = 
T_{m+n}(x) + \frac{1}{2^{2n}}T_{m-n}(x),\qquad m \geq n\neq 0 
\nonumber
\eea
Similarly,
\bea
U_m(x)T_n(x) = \left\{ \begin{array}{ll}
U_{m+n}(x) + \frac{1}{2^{2n}}U_{m-n}(x), \qquad m \geq n\neq 0 \\
U_{m+n}(x) - \frac{1}{2^{2m}}U_{n-m}(x), \qquad n \geq m\neq 0 
\end{array}
\right\}.
\nonumber
\eea
Analogous formulas to the ones above do not hold for all of the generalized Chebyshev polynomials. 
However, $P_{n\kp}$ play a similar role to the classical Chebyshev polynomials.

{\bf Lemma 1} $\quad$ {\it For any $m,n \in \BZ_+$ we have 
\bea
P_n(x)P_{m\kp}(x) = \left\{ \begin{array}{ll}
P_{n+m\kp}(x) + \frac{h_{m\kp}}{2}P_{n-m\kp}(x) \qquad n \geq m\kp\neq 0 \\
P_{n+m\kp}(x) + \frac{h_n}{2}P_{m\kp-n}(x) \qquad m\kp \geq n\neq 0
\end{array}
\right\}
\eea
}
Proof :

From section 1,
\bea
P_n(x) = e^{-n\chi_0}\frac{\vt(\up)}{\vt(\up-n\bh)}\left[e^{n\Omega(\pg_x)}\frac{\vt(\ux-n\bh)}{\vt(\ux)}+e^{-n\Omega(\pg_x)}\frac{\vt(\ux+n\bh)}{\vt(\ux)}\right],\nonumber
\eea
and,
\bea
P_{m\kp}(x) = e^{-m\kp \chi_0}\left[e^{m\kp \Omega(\pg_x)} + e^{-m\kp \Omega(\pg_x)}\right].\nonumber
\eea
Therefore,
\bea
P_n(x)P_{m\kp}(x) &=& e^{-(n+m\kp)\chi_0}\frac{\vt(\up)}{\vt(\up-n\bh)}\bigg[e^{(n+m\kp)\Omega(\pg_x)}\frac{\vt(\ux-n\bh)}{\vt(\ux)}+
e^{-(n+m\kp)\Omega(\pg_x)}\frac{\vt(\ux+n\bh)}{\vt(\ux)}\nonumber\\
&&\qquad +e^{(n-m\kp)\Omega(\pg_x)}\frac{\vt(\ux-n\bh)}{\vt(\ux)}+e^{-(n-m\kp)\Omega(\pg_x)}\frac{\vt(\ux+n\bh)}{\vt(\ux)}\bigg].\nonumber
\eea
From the periodicity of the theta functions we may add the vector $\pm m\kp \bh$ to any of the arguments of the theta functions without changing the equality. We make this change in the argument of the theta functions in such a way that we are left with (4.9).
$\Box$

Similarly, for the polynomials of the second kind, we have,
\bea
Q_n(x)P_{m\kp}(x) = \left\{ \begin{array}{ll}
Q_{n+m\kp}(x) + \frac{h_{m\kp}}{2}Q_{n-m\kp}(x) \qquad n \geq m\kp\neq 0 \\
Q_{n+m\kp}(x) - \frac{h_n}{2}Q_{m\kp-n}(x) \qquad m\kp \geq n\neq 0
\end{array}
.\right\}\nonumber
\eea
\noindent
Combining (4.9) and (4.7) we obtain,
\bea
P_{n\kp-j}(x) = \frac{1}{h_j}\left(P_{nk}(x)P_j(x)-\frac{Q_{nk}(x)Q_j(x)}{\psi^2(x)}\right) \qquad j = 1,..,\kp-1 .\nonumber
\eea
Similarly,
\bea
Q_{n\kp-j}(x) = \frac{1}{h_j}(Q_{n\kp}(x)P_j(x)-P_{n\kp}(x)Q_j(x)).\nonumber
\eea

\noindent
From (1.8) it follows from the periodicity properties of the theta functions that,
\bea
a_j = a_{\kp-j+1} \qquad j = 1,.., \kp-1. \nonumber
\eea
The $b_n$ also have a similar behaviour, and one can show that
\bea
b_j = b_{\kp-j+2} \qquad j= 1,.., \kp-1. \nonumber
\eea
However, this is not obvious from the form of the recurrence coefficients given in (1.9). To prove the above equality we start with the
three term recurrence relation (1.4) with $n=\kp-j+1$ and $x=1$, so that
\bea
P_{\kp-j+1}(1) = P_{\kp-j+2}(1)+b_{\kp-j+2}P_{\kp-j+1}(1)+a_{\kp-j+1}P_{\kp-j}(1). \nonumber
\eea  
From the expression for $P_{n\kp-j}(x)$ we have,
\bea
P_{n\kp-j}(1) = \frac{P_{n\kp}(1)P_j(1)}{h_j}.\nonumber
\eea
Substituting this into the three term recurrence relation we obtain,
\bea
P_{j-1}(1) = P_j(1) + b_{\kp-j+2}P_{j-1}(1) + a_{j-1}P_{j-2}(1) .\nonumber
\eea
Comparing this with (1.4) at $x=1$ and $n=j-1$ we see that $b_{\kp-j+2}=b_j$.

\begin{theorem}

Given a set $E$ such that the recurrence coefficients have period $\kp$, then 
the polynomial $P_\kp(x)/\Delta_\kp$ maps each sub-interval of $E$ onto 
$[-1,1]$.

\end{theorem}

{\it Proof:}
Take an arbitrary polynomial $M_\kp(x)$ of degree $\kp$ with leading coefficient $l_\kp$ where all of the zeroes lie in $(-1,1)$. Furthermore,  
\bea
|M_\kp(x_j)|\geq 1\nonumber\
\eea
where the $x_j$ are the stationary points, and, 
\bea
M_\kp(1)&=&1\nonumber\\
M_\kp(-1)&=&(-1)^\kp.\nonumber
\eea
Then the inverse image of $[-1,1]$ under $M_\kp$ is a set $J$ of at most $\kp$ disjoint intervals of the form $J:=[-1,\al]\cup_{j=1}^{\kp-2}[\bt_j,\al_{j+1}]\cup[\bt_{\kp-1},1]$, with $-1\leq \al_j \leq \bt_j \leq \al_{j+1} \leq 1$.
Let $\sigma_T(x)$ be the weight function associated with the Chebyshev 
polynomials,
\bea
\sigma_T(x) = \frac{1}{\pi}\frac{1}{\sqrt{1-x^2}}.\nonumber
\eea
 Following \cite{Ger:Van} we define a function $W(x)$ by, 
\bea
\frac{1}{\pi}\sqrt{\frac{\prod_{j=1}^{\kp-1}(x-\al_j)}{(1-x^2)\prod_{j=1}^{\kp-1}(x-\bt_j)}}=\sigma_A(x)=W(x)\sigma_T(M_\kp(x)).
\eea
Now, $M_\kp(x)$ has the property that it maps all of the branch points into $-1$ or $1$, and therefore,
\bea
M_\kp^2(x)-1 = l_\kp^2(x^2-1)\prod_{j=1}^{\kp-1}(x-\al_j)(x-\bt_j).\nonumber
\eea
Since $M_\kp(x)$ is a polynomial, the above equation imposes constraints  on 
the branch points. Substituting for $M_\kp^2(x)$ in (4.10) we have,
\bea
W(x) = l_\kp\prod_{j=1}^{\kp-1}(x-\al_j).\nonumber
\eea
Now
\bea
\frac{W(z)}{M_\kp(z)-x} = \frac{\prod_{i=1}^{\kp-1}(z-\al_i)}{\prod_{i=1}^\kp(z-M_i^{-1}(x))} = \sum_{i=1}^\kp\frac{w_i(x)}{z-M_i^{-1}(x)}
\eea
where the $M_i^{-1}(x)$ are a complete set of inverse branches of $M_\kp(x)$. Now,
\bea
W(z) &=& (M_\kp(z)-x)\sum_{i=1}^\kp\frac{w_i(x)}{z-M_i^{-1}(x)}\nonumber
\eea
and
\bea
W(M_{i}^{-1}(x))=w_i(x)\sum_{j \neq i}(M_i^{-1}(x)-M_j^{-1}(x)),\nonumber
\eea
and therefore,
\bea
w_i(x) = \frac{W(M_i^{-1}(x))}{M^{\pri}(M_i^{-1}(x))}.\nonumber
\eea
Also, by comparing coefficients of $z^{\kp-1}$ in (4.11),
\bea
\sum_{i=1}^\kp w_i(x) = 1.\nonumber
\eea   
From the orthogonality relationship of the Chebyshev polynomials we have, 
\bea
h_n^{(T)}\delta_{nm}&=&\int_{-1}^{1}T_n(z)T_m(z)\sigma_T(z)dz\nonumber\\
&=&\sum_{i=1}^\kp\int_{-1}^1T_n(z)T_m(z)\left(\frac{W(M_i^{-1}(z))}{M^{\pri}(M_i^{-1}(z))}\right)\sigma_T(z)dz\nonumber\\ 
&=&\int_JT_n(M_\kp(x))T_m(M_\kp(x))\frac{W(x)}{M^{\pri}_\kp(x)}\sigma_T(M_\kp(x))M_\kp^{\pri}(x)dx\nonumber\\
&=&\int_JT_n(M_\kp(x))T_m(M_\kp(x))\sigma_A(x)dx,\nonumber
\eea
where the last line follows from (4.10). Using the general method in \cite{Ger:Van} theorem 1, it can be shown that
\bea
\int_J x^jT_n(M_{\kp}(x))\sigma_A(x)dx = 0 \qquad j < n\kp. \nonumber
\eea
Therefore,
\bea
T_n(M_\kp(x)) = c_nP_{n\kp}(x) \qquad n \geq 0,\nonumber
\eea
where $c_n$ is some constant.
So for $n=1$,
\bea
M_\kp(x) = c_1P_\kp(x).\nonumber
\eea
Now, from (4.5),
\bea
P_\kp(1)=\Delta_\kp \nonumber
\eea
and therefore $c_1=\frac{1}{\Delta_\kp}$, and we have proven that $\frac{P_{\kp}(x)}{\Delta_\kp}$ maps $J$ into $[-1,1]$ $\kp$ -fold. Also, since for the polynomials orthogonal on the set $J$ we have the shown (4.5) to hold, we can conclude that the recurrence coefficients are periodic with period $\kp$ and therefore $J$ is in fact the set $E$ referred to in the statement of the theorem. 
$\Box$

\noindent
Remark 1. This theorem has appeared in a more general setting, see \cite{Peher1} and \cite{Peher2} where 
necessary and sufficient conditions were given for the existence of an extremal polynomial on a certain 
set $E$, which corresponds to our mapping polynomial. In general the orthogonal polynomials on these sets
have asymptotically periodic recurrence coefficients. These are found to be orthogonal with respect to 
a more general class of weight functions similar to $w(x)$ of the form,
\bea
\frac{\pi_{\kp-1}(x)}{y(x)}\nonumber
\eea
where $\pi_{\kp-1}(x)$ is a polynomial of degree $\kp-1$ chosen such that the weight function is 
positive on the interior of the set $E$.  In equation (1.12) of \cite{Peher2} the extremal polynomial of 
degree $n$, ${\cal T}_n(z)$, is $\frac{P_{\kp}(z)}{\Delta_{\kp}}$ where $n,$ the asymptotic period is 
identified with our $\kp$. 

As a consequence of this theorem, when the recurrence coefficients have period $\kp$ the set $E$ has
at most $\kp$ sub intervals. Recall that $g$ is the number of gaps in $[-1,1]$. Therefore, the 
maximum value of $g$ for this set is $\kp-1$.  

In figure 1 we show an example of a mapping polynomial for $\kp=3$.
\begin{figure}[h]
\includegraphics[0cm,0cm][2cm,6cm]{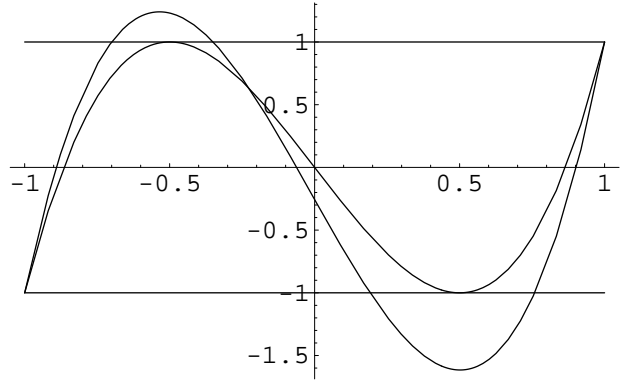}
\caption{General Form of $M_3(x)$ compared with the Chebyshev polynomial of degree $3$. The branch 
points are where the polynomial crosses the lines $y=-1$ and $y=1$. }
\end{figure}

We now consider some examples of sequences of orthogonal polynomials generated via the polynomial 
mappings $P_\kp(x)$. We choose small values of $\kp$ and investigate. $\kp=1$ corresponds to the 
Chebyshev polynomials as the recurrence coefficients are constant. 

For $\kp=2$, the mapping polynomial,
\bea
M_2(x) = \frac{1}{\Delta_2}P_2(x) = \frac{1}{\Delta_2}(x^2+ax+b_2),\nonumber
\eea
satisfies the following conditions
\bea
&&\frac{1}{\Delta_2}P_2(-1)=\frac{1}{\Delta_2}P_2(1)=1\nonumber\\
&&\frac{1}{\Delta_2}P_2(\al)=\frac{1}{\Delta_2}P_2(\bt)=-1.\nonumber
\eea
This implies that $M_2(x)$ has the following form 
\bea
M_2(x) = \frac{1}{\Delta_2}x^2+\left(1-\frac{1}{\Delta_2}\right),\nonumber
\eea 
with
\bea
\al^2&=&\bt^2\nonumber\\
\Delta_2&=&\frac{1-\al^2}{2} = \frac{1-\bt^2}{2}.\nonumber
\eea
Therefore,
\bea
M_2(x) = \frac{2}{1-\al^2}x^2+\left(1-\frac{2}{1-\al^2}\right),\nonumber
\eea
with the condition on the branch points, $\al=-\bt$.
The condition above, tells us which sets $E$ arise from a polynomial mapping of degree $2$. These are the same sets that give rise to periodic recurrence coefficients with period 2. To obtain the recurrence coefficients for these polynomials, we can substitute $\bt=-\al$ into the general genus 1 recurrence coefficients. The general expressions for the recurrence coefficients can be obtained in one of two ways. Appendix A shows how we can manipulate the theta function expressions for the recurrence coefficients into expressions involving the Jacobian elliptic functions. We are then able to invert the functions and in principle obtain the recurrence coefficients for any $n$ in terms of the branch points only. In this case we only need to find the first two since the others repeat. The second way to find the recurrence coefficients is to use the non-linear difference equations derived in section 2. Once we have the initial conditions we can iterate the non-linear difference equations to find the remaining coefficients. This method is more practical when we are considering higher periodicity. In any case, the recurrence coefficients are found to be, 
\bea
b_n&=&(-1)^n\al\nonumber\\
a_n&=&\frac{1}{4}(1-\al^2) \qquad n \geq 2.
\eea
With these coefficients, our polynomials become,
\bea
P_{2n}(x)=(\Delta_2)^nT_n\left(\frac{P_2(x)}{\Delta_2}\right)\nonumber
\eea
and
\bea
P_{2n+1}(x) = \frac{(x+\al)}{2}(\Delta_2)^n\left(T_n\left(\frac{P_2(x)}{\Delta_2}\right)+\frac{(x^2-1)}{\Delta_2}U_n\left(\frac{P_2(x)}{\Delta_2}\right)\right).\nonumber
\eea
Also the auxiliary polynomials become, 
\bea
&&S_1(x,n) = 2h_n(x+(-1)^{n+1}\al)\nonumber\\
&&G_2(x,n) = h_{n-1}(x+\al)(x-\al).\nonumber
\eea
From (3.10) we have,
\bea
f_1(x;n) = \frac{1}{2}\left(\frac{1}{x+(-1)^{n+1}\al}-\frac{1}{x-\al}\right).\nonumber
\eea
To evaluate $f_2(x;n)$ we must know the $c_0(n)$ and $\ga(n)$. For genus 1 polynomials we have from \cite{Chen-Lawrence},
\bea
c_0(n) = p_1(n)+n\frac{\al+\bt}{2}+\frac{\ga(n)-\al}{2}.\nonumber
\eea
where $p_1(n)$ is the coefficient of $x^{n-1}$ in $P_n(x)$.
In this case,
\bea
&&p_1(2n) = 0 \nonumber\\
&&p_1(2n-1) = \al. \nonumber
\eea
Therefore, 
\bea
c_0(n) = 0. \nonumber
\eea
Since $\ga(n)$ is the zero of $S_g(x,n)$,
\bea
\ga(n) = (-1)^n\al, \nonumber
\eea
and,
\bea
f_2(x;n) = \frac{nx}{x-\al}.\nonumber
\eea
For even $n$, the discriminant becomes ($n=2m$),
\bea
D[P_{2m}(x)] = (-1)^m2^{-(2m-1)(2m-2)-m}(2m)^{2m}(1-\al^2)^{m(2m-1)}T_{m}\left(\frac{\al^2+1}{\al^2-1} \right)\nonumber
\eea
and for odd $n$, ($n=2m+1$)
\bea
D[P_{2m+1}(x)] &=& (-1)^m 2^{-4m^2+m-1}(2m+1)^{2m}(1-\al^2)^{m(2m+1)}\nonumber\\
&&\qquad \qquad \qquad \left[T_m\left(\frac{\al^2+1}{\al^2-1}\right)-\frac{2}{1-\al^2}U_m\left(\frac{\al^2+1}{\al^2-1}\right)\right].\nonumber
\eea
The discriminant for the case when $n$ is odd cannot be obtained directly from (3.15) since the $P_{2n+1}(x)$ have zeroes at $x=-\al$. 
The problem is that of evaluating $P_{2n+1}^{\pri}(-\al)$. In general evaluating the derivative at one of the branch points is not 
straightforward. However, in this case it can be verified that, 
\bea
P_{2n+1}^{\pri}(-\al) = 2f_2(-\al;2n+1)Q_{2n+1}(-\al).\nonumber
\eea
whereas, for the zeroes that do not lie on the branch points we have,
\bea
P_{2n+1}^{\pri}(x_{j,2n+1}) = f_2(x_{j,2n+1};2n+1)Q_{2n+1}(x_{j,2n+1}) \nonumber
\eea
and the factor of $2$ does not appear. It is for this reason that the expression given in (3.15) remains valid in this case if we insert
a factor of $2$. 

The polynomials in this example can also be shown to satisfy the following differential equation,
\bea
P_n^{\pri\pri}-\frac{x^4-\al^2}{x(x^2-\al^2)(1-x^2)}P_n^{\pri}+\frac{n^2x^2}{(1-x^2)(x^2-\al^2)}P_n 
= 0 \qquad {\textrm {for even 
{\it n}}}, \nonumber
\eea
and,
\bea
P_n^{\pri\pri}-\frac{x^4-2\al x(1-x^2)+\al^2}{x(x^2-\al^2)(1-x^2)}P_n^{\pri}
+\frac{n^2x^4-2x\al+(n^2+1)x^3\al-\al^2}{x(1-x^2)(x-\al)(x+\al)^2}P_n 
= 0  \nonumber\\
\qquad \qquad \qquad \qquad{\textrm {for odd {\it n}}}.\nonumber
\eea

For the next example we choose $\kp=3$. Therefore our polynomials can be 
described in terms of a polynomial mapping of degree 3. The highest genus 
permitted in this case is 2. This corresponds to three disjoint intervals. 
Genus one cases arise if we close one of the gaps. In the last example we saw 
that the only sets that gave rise to a period of $2$ were symmetric around 
the origin, that  is $\al=-\bt$. The aim of this example is to determine which sets 
give rise to recurrence coefficients with period 3. We determine the 
recurrence coefficients, the $\ga$ points and under a certain condition, the 
$c_j(n)$.

We have a mapping polynomial $M_3(x)$ with leading coefficient $\frac{1}{\Delta_3}$ that satisfies the following equations,
\bea
&&M_3(\pm 1) = \pm 1\nonumber\\
&&M_3(\al_1) = M_3(\bt_1) = 1\nonumber\\
&&M_3(\al_2) = M_3(\bt_2) = -1.\nonumber
\eea
Solving these six equations we have,
\bea
M_3(x) = \frac{x^3}{\Delta_3}+ax^2+\left(1-\frac{1}{\Delta_3}\right)x-a,\nonumber
\eea
with,
\bea
a = \frac{-(\bt_1+\al_1+1)}{(\bt_1+1)(\al_1+1)} = \frac{-(\al_2+\bt_2-1)}{(\bt_2-1)(\al_2-1)}, 
\eea
and,
\bea
\Delta_3 = (\bt_1+1)(\al_1+1) = (\bt_2-1)(\al_2-1).
\eea

Since we can characterize the mapping polynomial in terms of two branch points, we have a system with two unknowns. We describe the system in terms of $\al_1$ and $\bt_1$. 
The form of the mapping polynomial restricts $\al_1$ and $\bt_1$ to certain intervals on $[-1,1]$. If we eliminate $\al_2$ from (4.13) and (4.14) we get,
\bea
\bt_2 = \frac{\al_1+\bt_1+2 \pm \sqrt{(\bt_1-\al_1)^2-4(\al_1+\bt_1+1)}}{2}\nonumber
\eea
hence,
\bea
(\al_1-\bt_1)^2-4(\al_1+\bt_1+1) \geq 0. 
\eea
Therefore we have,
\bea
\al_1 \leq \bt_1+2-4\sqrt{\frac{\bt_1+1}{2}}.\nonumber
\eea 
This inequality, coupled with the inequality,
\bea
\al_1 \leq \bt_1 \nonumber
\eea
implies that $\al_1$ lies on the interval $[-1,-\frac{1}{2}]$. Therefore $\bt_1$ lies on the interval $[\al_1,\bt_{1,max}]$ where,
\bea
\bt_{1,max} = 2+\al_1-4\sqrt{\frac{\al_1+1}{2}}.\nonumber
\eea
The value of $\bt_{1,max}$ follows from (4.15) with the equality sign.
With these conditions on $\al_1$ and $\bt_1$ we have,
\bea
\al_2 = 1+\frac{\al_1+\bt_1}{2}-\frac{1}{2}\sqrt{(\bt_1-\al_1)^2-4(\al_1+\bt_1+1)}
\eea
and,
\bea
\bt_2 = 1+\frac{\al_1+\bt_1}{2}+\frac{1}{2}\sqrt{(\bt_1-\al_1)^2-4(\al_1+\bt_1+1)}.
\eea
For $g=2$ we evaluate the first three recurrence coefficients in the most 
general case in terms of $\al_1$, $\al_2$, $\bt_1$ and $\bt_2$, using (1.13) 
and (1.14). Then we substitute for $\al_2$ and $\bt_2$ using (4.16) and (4.17)
and find,
\bea
a_1 &=& \frac{1}{2}(1+\al_1)\left(2-\al_1+\bt_1+\sqrt{(\bt_1-\al_1)^2-4(1+\al_1+\bt_1)}\right)\nonumber\\
a_2 &=& \frac{1}{16}\left(\al_1-\bt_1-2+\sqrt{(\bt_1-\al_1)^2-4(1+\al_1+\bt_1)}\right)^2\nonumber\\
a_3 &=& \frac{1}{4}(1+\al_1)\left(2-\al_1+\bt_1+\sqrt{(\bt_1-\al_1)^2-4(1+\al_1+\bt_1)}\right),
\eea
and,
\bea
b_1 &=& \frac{1}{2}\left(\bt_1-\al_1+\sqrt{(\bt_1-\al_1)^2-4(1+\al_1+\bt_1)}\right)\nonumber\\
b_2 &=& \frac{1}{4}\left(2+3\al_1+\bt_1-\sqrt{(\bt_1-\al_1)^2-4(1+\al_1+\bt_1)}\right)\nonumber\\
b_3 &=& \frac{1}{4}\left(2+3\al_1+\bt_1-\sqrt{(\bt_1-\al_1)^2-4(1+\al_1+\bt_1)}\right).
\eea
We also find that, for $n \geq 2$,
\bea
\gamma_{1,2}(n) &=& \frac{1}{2}\Bigg(1+\al_1+\bt_1-b_{n+1}\nonumber\\
&\mp&\sqrt{1-2(\al_1+\bt_1)+(\al_1-\bt_1)^2-4(a_n+a_{n+1})+2b_{n+1}(1+\al_1+\bt_1-\frac{3}{2}b_{n+1})}\Bigg). \nonumber 
\eea
Expressions for the polynomials are easily obtained from our equations in the previous sections. Expressions for the discriminant and the differential equation can in principle be obtained using the constants that we have derived in this section.

A highly symmetric case corresponding to $g=2$  arises when 
$\bt_1+\al_1 = -1$.
With this, the $a$ defined by (4.13) is $0$, which implies that
\bea
\al_2+\bt_2=1.\nonumber
\eea 
Substituting these into (4.14), gives,
\bea
\bt_2=1+\bt_1.\nonumber
\eea
This case corresponds to the first gap being symmetric about the point $-\frac{1}{2}$ and the second gap being symmetric about the point $\frac{1}{2}$. 
In this case the square roots appearing in the recurrence coefficients vanish leaving (with $\al$ = $\al_1$),
\bea
a_1 &=& -2\al(1+\al)\nonumber\\
a_2 &=& \frac{1}{4}\nonumber\\
a_3 &=& -\al(1+\al)\nonumber
\eea
and
\bea
b_1 &=& -(1+2\al)\nonumber\\
b_2 &=& \frac{1}{2}+\al\nonumber\\
b_3 &=& \frac{1}{2}+\al\nonumber
\eea
Also,
\bea
G_3(x,n) = h_{n-1}(x^3+(2a_n-1-\al-\al^2)x+2a_n(b_n+b_{n+1}))\nonumber
\eea
and
\bea
S_2(x,n) = 2h_n(x^2+b_{n+1}x+a_n+a_{n+1}+b_{n+1}^2-1-\al-\al^2).\nonumber
\eea
This gives
\bea
\ga_{1,2}(n) = \frac{-b_{n+1}\mp\sqrt{b_{n+1}^2-4(a_n+a_{n+1}+b_{n+1}^2-1-\al-\al^2)}}{2}. \nonumber
\eea
The first three polynomials are,
\bea
P_1(x) &=& x+2\al+1\nonumber\\
P_2(x) &=& x^2+(\al+\frac{1}{2})x-\frac{1}{2}\nonumber\\
P_3(x) &=& x^3-(1+\al+\al^2)x.\nonumber
\eea
From (3.9), (3.10) and (3.12) we can calculate $f_1(x;n)$ and $f_2(x;n)$. 


Returning to the previous non-symmetric genus 2 example, we add the restriction that one of the gaps is closed. We take the case where $\al_2=\bt_2$. In this
situation (4.13) and (4.14) simplify to 
\bea
a = \frac{-(\bt_1+\al_1+1)}{(\bt_1+1)(\al_1+1)} = \frac{-(2\al_2-1)}{(\al_2-1)^2}, \nonumber
\eea
and
\bea
\Delta_3 = (\bt_1+1)(\al_1+1) = (\al_2-1)^2.\nonumber
\eea
Therefore, eliminating $\al_2$ from the above we obtain,
\bea
(\bt_1+1)(\al_1+1) = \left(\frac{\bt_1+\al_1}{2}\right)^2,\nonumber
\eea
with the solutions,
\bea
\bt_1 = \al_1+2\pm 4\sqrt{\frac{\al_1+1}{2}}.\nonumber
\eea
Since $\al_1 \leq \bt_1$ we must have
\bea
\bt_1 = \al_1+2-4\sqrt{\frac{\al_1+1}{2}}.
\eea 
This relation governs the variation of $\bt_1$ as $\al_1$ varies over the 
interval $[-1,-\frac{1}{2}]$. It describes a contour in the $(\al_1,\bt_1)$ 
plane. Figure (2) shows a graphical representation of this contour with 
similar curves for higher values of $\kp$.  
Had we have chosen the condition $\al_1=\bt_1$ for the other $g=1$ case we 
would have obtained,
\bea
\bt_2=-2+\al_2+4\sqrt{\frac{1-\al_2}{2}},
\eea 
where $\al_2$ varies over $[-1,\frac{1}{2}]$.
Using the method outlined in the appendix we are able to compute the recurrence coefficients for $g=1$ explicitly. From now on we refer to $\al_1$ and $\bt_1$ as just $\al$ and $\bt$ respectively. Upon substituting for $\bt$ given by (4.21) into the general recurrence coefficients we obtain the following, 
\bea
a_1 &=& 2(\al+1)\left(1-\sqrt{\frac{1+\al}{2}}\right)\nonumber\\
a_2 &=& \frac{\al+3}{2}-2\sqrt{\frac{\al+1}{2}}\nonumber\\
a_3 &=& (\al+1)\left(1-\sqrt{\frac{1+\al}{2}}\right)\nonumber
\eea
and,
\bea
b_1 &=& 1-2\sqrt{\frac{1+\al}{2}}\nonumber\\
b_2 &=& (\al+1)-\sqrt{\frac{1+\al}{2}}\nonumber\\
b_3 &=& (\al+1)-\sqrt{\frac{1+\al}{2}}.\nonumber
\eea
The recurrence coefficients repeat in blocks of three however we must remember that $a_{3n+1} = \frac{a_1}{2}, \;\; n \geq 1$. Furthermore,
\bea
c_0(3n) &=& n\left(\al-2\sqrt{\frac{1+\al}{2}}\right)\nonumber\\
c_0(3n+1) &=& \frac{1}{2}\left((1+2n)\al-(1+4n)\sqrt{\frac{1+\al}{2}}\right)\nonumber\\
c_0(3n+2) &=& \frac{1}{2}\left((1+2n)\al-(3+4n)\sqrt{\frac{1+\al}{2}}\right),\nonumber
\eea
and,
\bea
\ga(3n) &=& \al\nonumber\\
\ga(3n+1) &=& -\sqrt{\frac{1+\al}{2}}\nonumber\\
\ga(3n+2) &=& -\sqrt{\frac{1+\al}{2}}.\nonumber
\eea
These constants, listed above, allow one to calculate the differential equations and the discriminant given in the previous section. 

{\bf General $\kp$}

For notational convenience we now label as $E_{\kp}$ the sets $E$ that give rise to period $\kp$ 
recurrence coefficients. Equivalently these are the sets for which the minimum degree of polynomial 
mapping from $E_{\kp}$ onto $[-1,1]$ is $\kp$. From our earlier examples,
\bea
E_{1} &=& [-1,1] \nonumber\\
E_{2} &=& [-1,\al] \cup [-\al,1]  \qquad \al \in [-1,0]\nonumber\\
E_{3} &=& [-1,\al] \cup [\bt,u(\al,\bt)] \cup [v(\al,\bt),1] \qquad \al \in [-1,-\frac{1}{2}] \;\;\; \bt \in [\al,\bt_{max}],\nonumber  
\eea
where,
\bea
u(\al,\bt) &=&  1+\frac{\al+\bt}{2}-\frac{1}{2}\sqrt{(\bt-\al)^2-4(\al+\bt+1)} \nonumber\\
v(\al,\bt) &=&  1+\frac{\al+\bt}{2}+\frac{1}{2}\sqrt{(\bt-\al)^2-4(\al+\bt+1)}.\nonumber
\eea
Each set $E_{\kp}$ is parametrized in terms of $\kp-1$ variables. As $\kp$ 
increases it becomes intractable to solve explicitly for the relationships 
between the branch points. However, the branch points are constrained on 
certain curves in $[-1,1]^{2g}$ where the co-ordinate axes of the $2g$ 
cube label the position of the branch points. 

In what follows, we describe the constraints on the mapping polynomial. 
This will provide a set of equations that define the curves of 
constant periodicity in the $2g$ cube. We also show how to determine the 
smallest $\kp$ for a given configuration of the branch points.

For general $\kp$, the $2\kp$ conditions imposed on the mapping polynomial,
\bea
M_{\kp}(x) = a_{\kp}x^{\kp} + a_{\kp-1}x^{\kp-1} + ... + a_1 x + a_0, 
\nonumber
\eea
are 
\bea
&& M_{\kp}(-1) = (-1)^{\kp}\nonumber\\
&& M_{\kp}(1) = 1 \nonumber\\
&& M_{\kp}(\al_j) = M_{\kp}(\bt_j)  = (-1)^{\kp+j} \qquad j = 1,..., \kp-1. 
\eea
Therefore $\kp+1$ of these will determine the coefficients in terms of $\kp+1$ 
branch points. Using the first two equations of (4.22) and any $\kp-1$ of the 
remaining ones we can solve for the coefficients of the mapping polynomial
in terms of $1$, $-1$ and $\kp-1$ free parameters. The other $\kp-1$ branch
points are related to the free parameters by the remaining equations of 
(4.22).   
These equations determine a surface in the $2g$ cube. The surface obtained in 
this way will correspond to maximum $g=\kp-1$. Closing one of the gaps, that  is setting  
$\al_j=\bt_j$ for some $j \leq \kp-1$, will determine surfaces corresponding to that value of 
$\kp$ for smaller values of $g$. 

If we close all of the gaps by setting $\al_j = \bt_j$ for all $j=1...\kp-1$, then the mapping
polynomial becomes,
\bea
M_{\kp}(x) = 2^{\kp-1}T_{\kp}(x).\nonumber
\eea
If we denote by $\tilde{x}_j$ the point where $\al_j=\bt_j$, then we must have 
\bea
2^{\kp-1}T_{\kp}(\tilde{x}_j) = \pm 1\nonumber
\eea
from the properties of the mapping polynomial. Since $T_n(x) = 2^{1-n}\cos(n\theta), \;\; x= \cos
(\theta)$ we see that the set of points $\{\tilde{x}_j\}_{j=1}^{\kp-1}$ are the $\kp-1$ stationary
points of the Chebyshev polynomial of degree $\kp$. Therefore, each gap i.e. 
$(\al_i,\bt_i)_{i=1}^{\kp-1}$ can be said to be 'centred' (not necessarily symmetrically) around 
one of the stationary points of the Chebyshev polynomial of degree $\kp$. The most simple example 
of this is the $\kp=2$ case where the gap is centred around the point $x=0$. In the $\kp=3$, $g=2$
case, the two gaps are centred around $x=\frac{1}{2}$ and $x=-\frac{1}{2}$. The stationary points
of $T_4(x)$ are $x=-\frac{1}{\sqrt{2}}$, $x=0$, and $x=\frac{1}{\sqrt{2}}$. So for a polynomial 
mapping of degree $4$ we have three distinct $g=1$ cases corresponding to one gap which is opened
 around one of the stationary points. 


Consider now the $\kp=2$ mapping polynomial, $\frac{P_2(x)}{\Delta_2}$, which maps 
$[-1,\al] \cup [-\al,1]$ to $[-1,1].$
From (4.4),
\bea
P_4(x) = \Delta_2^2T_2\left(\frac{P_2(x)}{\Delta_2}\right).\nonumber
\eea
Since, $\Delta_2^2 = 2\Delta_4,$ and from the fact that the classical Chebyshev polynomials map
$[-1,1]$ to $[-1,1]$ we see that, $\frac{P_4(x)}{\Delta_4}$ also maps $[-1,\al] \cup [-\al,1]$ to 
$[-1,1]$. Therefore, it is this configuration of the branch points which corresponds to the 
$\kp=4$, $g=1$ case where the gap is centered around $x=0$. Indeed for any mapping polynomial of 
even degree, one of the $g=1$ cases corresponds to this configuration, since a Chebyshev 
polynomial of even degree has a stationary point at zero. 

The same reasoning as above applies for general $\kp$. that
is, if 
$\frac{P_{\kp}(x)}{\Delta_{\kp}}$ is a mapping polynomial for a certain configuration of the 
branch points, then so are the polynomials $\frac{P_{n\kp}(x)}{\Delta_{n\kp}}$. Therefore
amongst all of the configurations that correspond to a mapping of degree $n\kp$, there can be 
found those that correspond to degree $j\kp$, where $j$ divides $n$.

If $\tilde{x}_l$ is a stationary point of $\widehat{T}_\kp(x):= 2^{\kp-1}T_{\kp}(x)$, then 
\bea
|\widehat{T}_\kp(\tilde{x}_l)| = 1\nonumber
\eea
and 
\bea
|\widehat{T}_{j\kp}(\tilde{x}_l)| = 1, \nonumber
\eea
which indicates that $\tilde{x}_l$ is also a stationary point of $\widehat{T}_{j\kp}(x)$. For 
example $\widehat{T}_{12}(x)$ has amongst 
its stationary points, the stationary points of 
$\widehat{T}_6(x)$, $\widehat{T}_4(x)$, $\widehat{T}_3(x)$ 
and $\widehat{T}_2(x)$. 

So as an example, a mapping polynomial of degree $12$ 
will give rise to eleven $g=1$, single 
parameter cases of the generalized Chebyshev polynomials. 
Six of these will already be the same 
as those that arise from a mapping of degree $6$ and 
four will be the same as those that arise 
from a mapping of degree $4$. All of these cases cover 
those that arise from mappings of degree 
$2$ and $3$. Therefore a mapping of degree $12$ will 
give rise to two $g=1$ cases that cannot
be generated from polynomial mappings of lower degree. 

The upshot of the above discussion is that we can 
uniquely label each surface in the $2g$ cube 
by a set of points corresponding to the stationary 
points of the Chebyshev polynomials of 
any degree. We do this in the following way.  

If we denote by $s_l$ the ordered set of all 
stationary points of the Chebyshev polynomial of 
degree l, then we can label a surface by taking any 
subset of the $s_l$. The number of elements
in a subset is the number of gaps $g$ in the interval 
$[-1,1]$. For example,
\bea
s_2 &=& \left\{ 0\right\}\nonumber\\
s_3 &=& \left\{ -\frac{1}{2}, \frac{1}{2}\right\}
\nonumber\\
s_4 &=& \left\{-\frac{1}{\sqrt{2}},0,\frac{1}
{\sqrt{2}}\right\}
\nonumber\\
s_5 &=& \left\{-\frac{1}{4}(\sqrt{5}+1),-\frac{1}{4}
(\sqrt{5}-1), \frac{1}{4}(\sqrt{5}-1),
\frac{1}{4}(\sqrt{5}+1) \right\}. 
\nonumber
\eea
Therefore, the set,
\bea
\left\{-\frac{1}{\sqrt{2}},\frac{1}{\sqrt{2}}\right\}\nonumber
\eea
corresponds to the $\kp=4$ $g=2$ surface where both gaps are centred around the two points 
$-\frac{1}{\sqrt{2}}$ and $\frac{1}{\sqrt{2}}$. This is a two parameter curve. Those curves  
shown in figure 2 can be labelled by one number which is a stationary point of a Chebyshev 
polynomial. It is clear from the discussion so far that the number ${\cal N}$ given by,
\bea
{\cal N} = {\kp-1 \choose g},\nonumber
\eea 
is the number of surfaces in the $2g$ cube associated with the integer $\kp$.

So far, we have discussed in general terms how to describe the surfaces 
corresponding to a particular value of $\kp$. On the other hand, the problem of finding $\kp$ 
from a given set $E$ appears to be a very difficult numerical one. The condition for periodicity 
of the recurrence coefficients is 
\bea
\kp \bh \in \BZ^g\setminus \{0^g\} .
\eea 
where,
\bea
\bh_j = \frac{1}{2\pi i}\int_{b_j}d\Omega, \qquad j=1,...,g.\nonumber
\eea
If $-\bh_j$ are rational then $\kp$ is the least 
common multiple of the denominators
of these national numbers. For example, if
\bea
\bh = \left( \begin{array}{lll} 
-\frac{1}{3} \\
-\frac{1}{5} \\
-\frac{1}{6}
\end{array}
\right)
\nonumber
\eea
then $\kp = 30$. However, if at 
least one of the components of the vector $\bh$ is an
irrational  number then it is impossible to find an integer
$\kp$ that satisfies condition (4.25). These  correspond to
the cases where the recurrence coefficients never repeat
themselves or equivalently, the sets for which there is no
polynomial mapping from this set onto $[-1,1]$. 

Note that in general the values of $\bh_j$ are 
subject to certain conditions. In 
\cite{Chen-Lawrence}, section 10, a physical 
interpretation was given; the values $-\bh_j$ were 
seen to be the proportion of the overall charge 
lying on 

the interval $[\bt_{j-1},\al_j]$, where $\bt_0:=-1$, in the equilibrium charge density for the 
set $E$, when the total charge was normalized to 1, 
that  is,
\bea
\int_E \sigma(x)dx = 1. \nonumber
\eea 
Since the total charge on all of the intervals must 
add up to 1 the following inequality is 
satisfied by the components of $\bh$,
\bea
0 < \sum_{j=1}^g -\bh_j < 1 .\nonumber
\eea

\begin{figure}[h]
\includegraphics{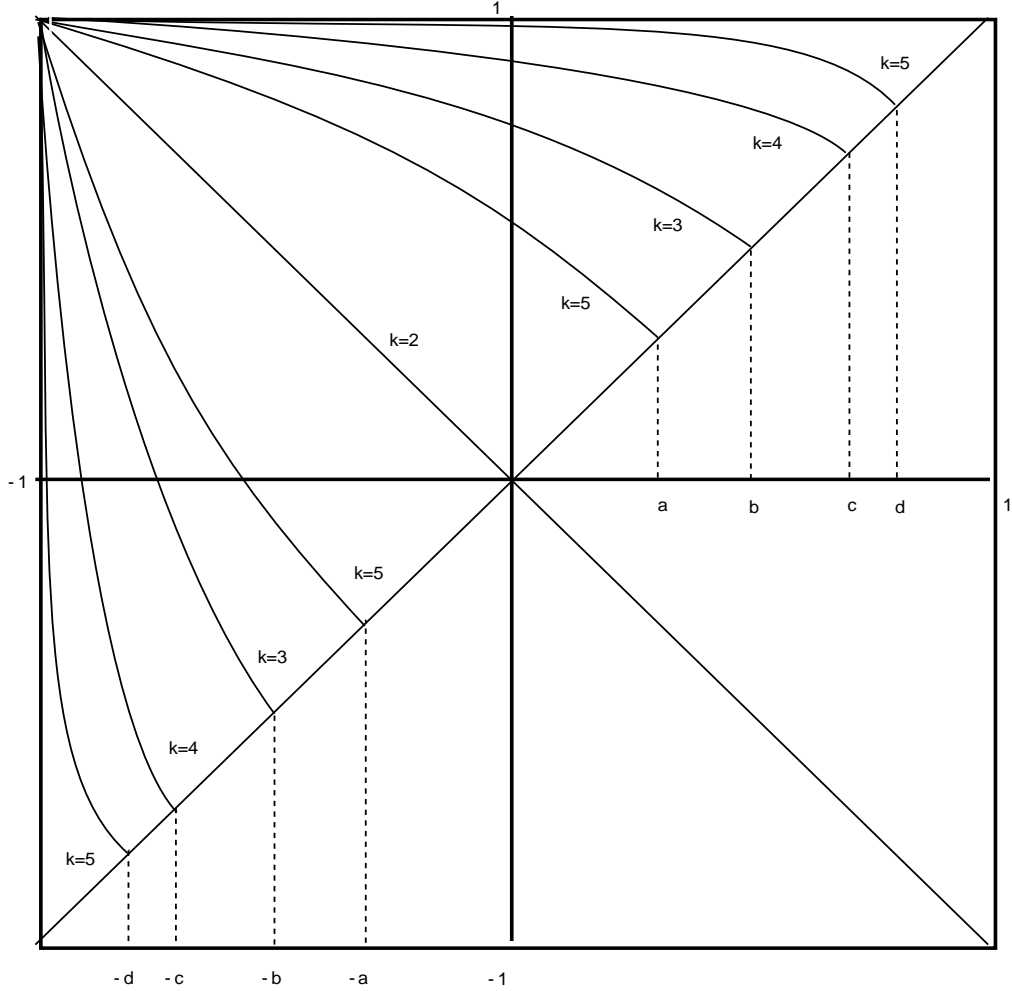}
\caption{Contours of constant periodicity in 
the $(\al,\bt)$ plane for genus 1. Here $a=(-1+{\sqrt {5}})/4,$ $b=1/2,$ $c=1/{\sqrt {2}}$ 
and $d=(1+{\sqrt {5}})/4.$ This is not a
precise plot, it is representative of the relative position of the contours.}
\end{figure}

To finish this section we refer to figure $2$. Recall that each point in the plane represents a 
configuration of the branch points for $g=1$ and therefore a sequence of orthogonal polynomials. 
Using our notation, 
if $\{\tilde{x}_l\}$ represents a curve in this plane with integer $\kp$ then $\{-x_i\}$ 
represents a different curve but with the same integer $\kp$. We ask the question ``Is there a 
relationship between the polynomials associated with the curves $\{x_i\}$ and $\{-x_i\}$?'' 
A similar question can also be asked for $g > 1$.
The relation between these cases is given in the most general form by theorem $1.2$. However, 
since the recurrence coefficients are periodic, the elements appearing in the determinant 
simplify:
\bea
P_{n\kp+j}(\bt_l) = \frac{1}{2}P_{n\kp}(\bt_l)P_j(\bt_l) 
\qquad 1 \leq j \leq \kp-1\nonumber
\eea
and
\bea
Q_{n\kp+j}(\al_l) = \frac{1}{2}P_{n\kp}(\al_l)Q_j(\al_l) 
\qquad 1 \leq j \leq \kp-1\nonumber
\eea
where,
\bea
P_{n\kp}(\al_l) = P_{n\kp}(\bt_l) = 
\frac{\Delta_\kp^n(-1)^{n(l+\kp)}}{2^{n-1}}.
\nonumber
\eea 
As an example we show the relation between the 
polynomials associated with the curves 
$\{-\frac{1}{2}\}$ and $\{\frac{1}{2}\}$ which are 
the two $\kp=3$ $g=1$ curves. Using the 
notation of theorem 1.2, we denote the polynomials 
associated with the set 
$\left\{-\frac{1}{2}\right\}$ as $P_n(x)$, 
the polynomials associated with the set 
$\left\{\frac{1}{2}\right\}$ as $\tilde{P}_n(x)$, 
and the respective mapping polynomials
 as $M_3(x)$ and $\tilde{M}_3(x)$. It is easy to verify that,
\bea
M_3(x) = (-1)\tilde{M}_3(-x).\nonumber
\eea
Consequently,
\bea
\tilde{P}_{3n}(x) = (-1)^{3n}P_{3n}(-x).\nonumber
\eea
However, for the intermediate polynomials we must use the expressions given in theorem 1.2, 
which give,
\bea
\tilde{P}_{3n+1}(x) &=& \frac{(-1)^{3n}}{\left(x+\al+2-4\sqrt{\frac{\al+1}{2}}\right)}\Bigg[P_{3n+2}(-x)-\left(\sqrt{\frac{\al+1}{2}}-1\right)P_{3n+1}(-x)\nonumber\\
&& \qquad \qquad \qquad \qquad \qquad \;\;\;+\left((\al+3)\sqrt{\frac{\al+1}{2}}-2(\al+1)\right)P_{3n}(-x)\Bigg],\nonumber
\eea
and,
\bea
\tilde{P}_{3n+2}(x) &=& \frac{(-1)^{3n+1}}{\left(x+\al+2-4\sqrt{\frac{\al+1}{2}}\right)}\Bigg[P_{3(n+1)}(-x)+\left(\sqrt{\frac{\al+1}{2}}\right)P_{3n+2}(-x)\nonumber\\
&& \qquad \qquad \qquad \qquad \qquad \qquad \;\; -\frac{1}{2}\left(1+\al-2\sqrt{\frac{1+\al}{2}}\right)P_{3n+1}(-x)\Bigg].\nonumber
\eea

{\bf 4b. General form of the polynomials}

In this section we show some plots of the polynomials for small $\kp$ and prove a theorem 
regarding the zeroes. The Chebyshev polynomials $\widehat{T}_n(x)$ is bounded between $1$ and 
$-1$. These are the mapping polynomials from $[-1,1]$ to $[-1,1]$. 
As a consequence of (4.5) we see that $|\widehat{P}_{n\kp}(x)| \leq 1$ for $x \in E$, 
where,
\bea
\widehat{P}_{n\kp}(x)=\frac{2^{n-1}}{(\Delta_\kp)^n}P_{n\kp}(x).\nonumber
\eea
Outside of the set we no longer have this bound. Therefore we would expect the polynomials to 
appear graphically as in figure 3. In figures 4, 5 and 6  we plot some of the intermediate 
polynomials $\widehat{P}_{n\kp+j}(x)$ for various values of $j$ and $\kp$, where,
\bea
\widehat{P}_{n\kp+j}(x)=\frac{2^{n}}{(\Delta_\kp)^n\sqrt{2h_j}}P_{n\kp+j}(x).\nonumber
\eea
In $E$ they are bounded above and below by some function. It is clear from the representation 
(1.18) that the form of the envelope is $\pm \wh{\rho}_{n\kp+j}(x)$, for $x \in E$ where, 
\bea
\widehat{\rho}_{n\kp+j}(x) &=& \frac{2^{n}}{(\Delta_\kp)^n\sqrt{2h_j}}\sqrt{\frac{S_g(x,n\kp+j)}{\prod_{l=1}^g(x-\al_l)}},\nonumber\\
&=& \sqrt{\frac{\prod_{j=1}^g(x-\ga_j(n\kp+j))}{\prod_{j=1}^g(x-\al_j)}}.\nonumber
\eea
Note that the quantity under the square root is strictly positive for $x \in E$.

\begin{figure}[htb]
\centering
\includegraphics{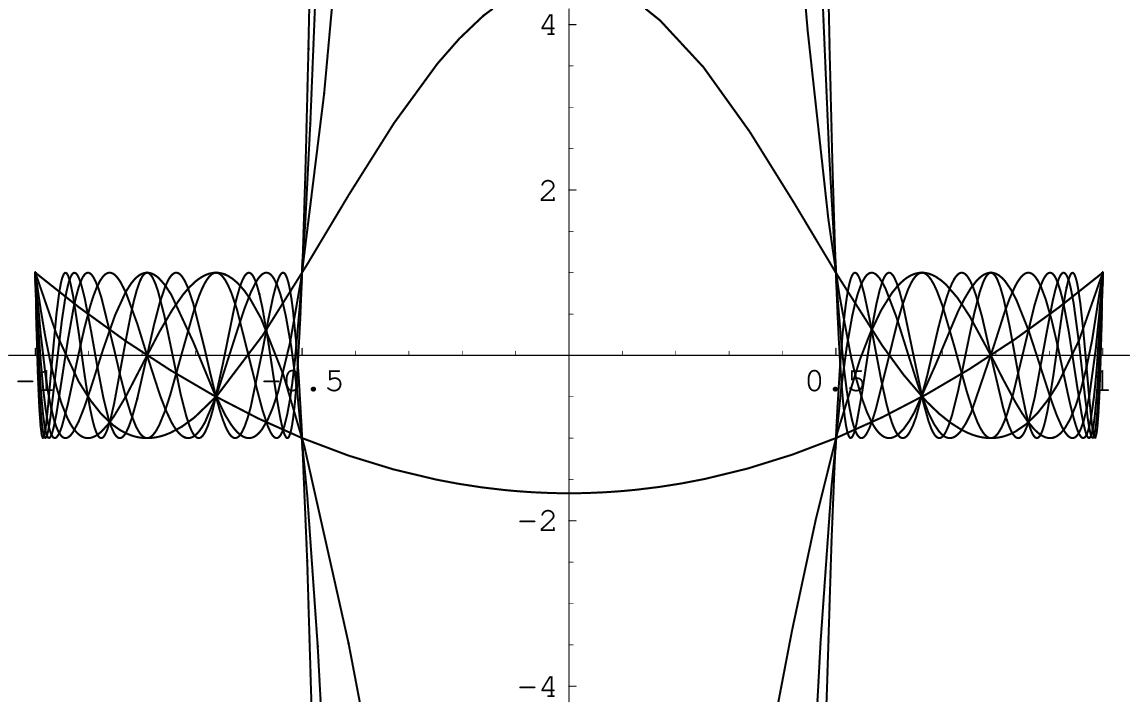}
\caption{$\kp$=2, j=0  n=1..8, $\al=-0.5$}
\end{figure}

\begin{figure}[htb]
\centering
\includegraphics{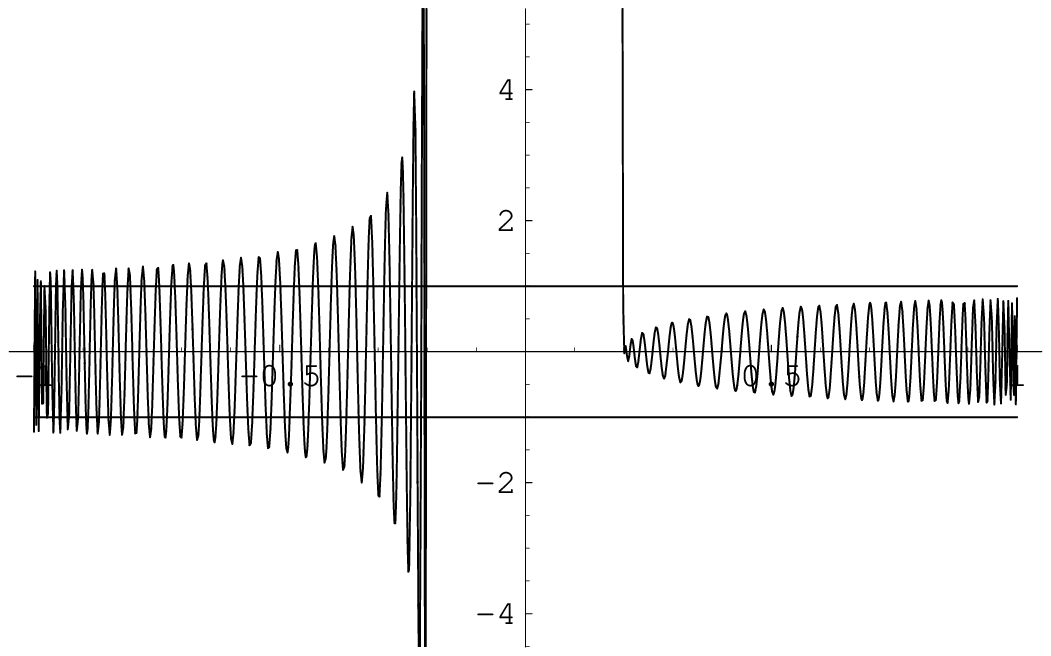}
\caption{$\kp$=2, j=1, n=65, $\al=-0.2$}
\end{figure}

\begin{figure}[htb]
\centering
\includegraphics{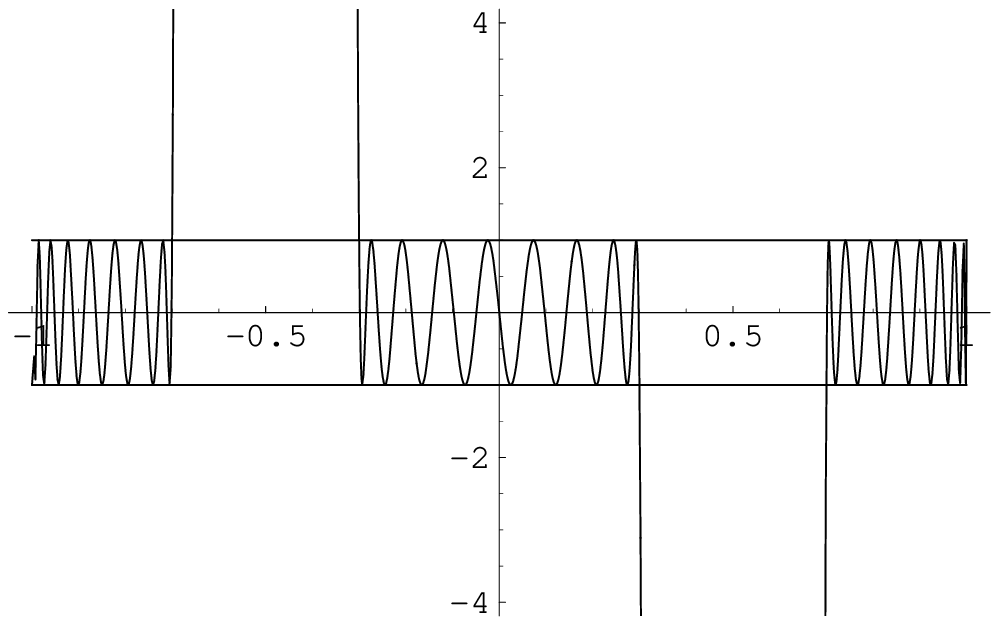}
\caption{The symmetric case for 
$\kp$=3, j=0, n=17, $\al=-0.7$}
\end{figure}

\begin{figure}[htb]
\centering
\includegraphics{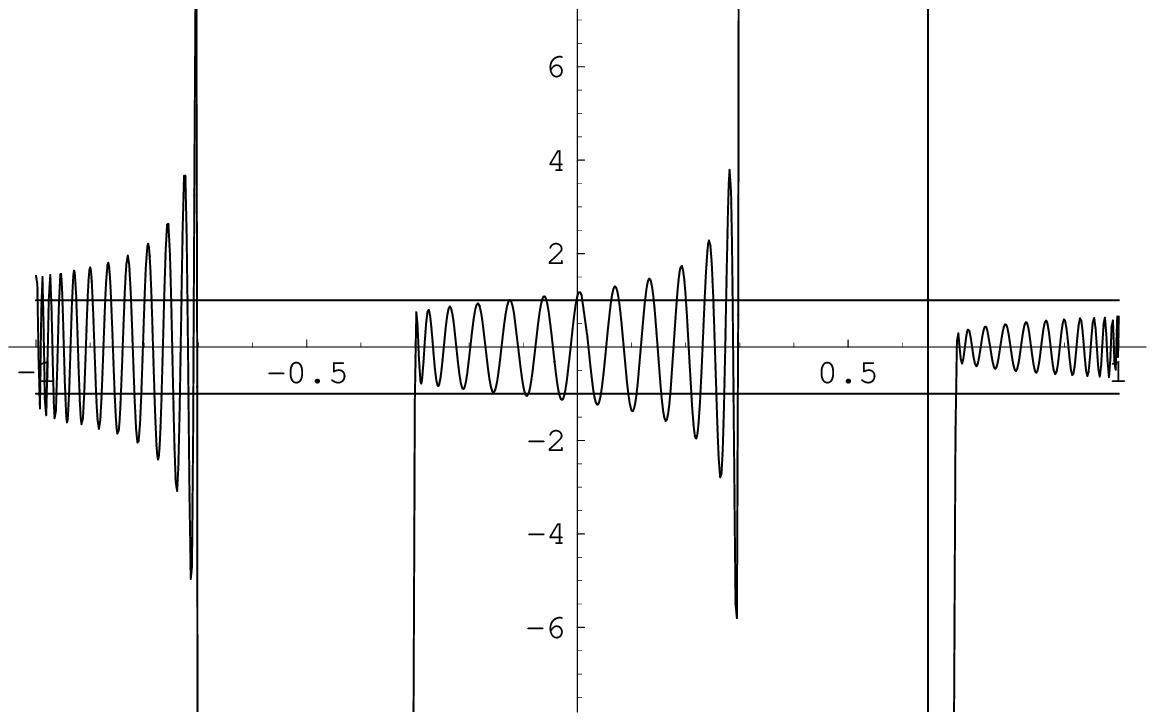}
\caption{The symmetric case for 
$\kp$=3, j=1, n=25, $\al=-0.7$}
\end{figure}

We already know something about the zeros of 
the generalized chebyshev polynomials in the most 
general case from \cite{Chen-Lawrence}. In the 
polynomial mapping case we can say a little more 
about the zeros. 

\begin{theorem}

All of the zeros of $P_{n\kp}(x)$ lie inside $E$, 
and in each of the $\kp$ sub-intervals that 
make up $E$ lie $n$ zeros.

\end{theorem}

{\it Proof:}
All of the zeros of the Chebyshev polynomials are 
contained in $(-1,1)$. From (4.5), and since 
$\frac{P_\kp(x)}{\Delta_\kp} \in [-1,1]$ if and only 
if $x \in E$, we see that  all of the zeros 
of $P_{n\kp}(x)$ lie inside $E$ and each subinterval $E_i$ 
contains $n$ zeros.  

The zeros of $Q_{n\kp}(x)$ behave similarly except 
in this case there are $n-1$ zeros in each 
interval $E_i$ and there are also zeros at 
every $\al_j$, $j=1,..,\kp-1$. This is easily verified 
by checking the form of (4.6). From (4.6) the $n-1$ zeros 
in each $E_i$ come for the zeros of 
$U_n\left(\frac{P_\kp}{\Delta_\kp}\right)$. In each 
interval $E_i$, these zeros interlace with 
the zeros of $P_{n\kp}(x)$.

\begin{theorem}

Let $\{x_{i,n}\}$ denote the zeros of $P_{n}(x)$ 
and $\{y_{i,n}\}$ denote the zeros of $Q_n(x)$.
Consider the following two sets,
\bea
A&:=& \{\{x_{i,n\kp}\},\{x_{i,j}\}\}\nonumber\\
B&:=& \{-1,1,\bt_1....\bt_{\kp-1},\{y_{i,j}\},
\{y_{i,n\kp}\}\}\setminus\{\al_1...\al_{\kp-1}\}.
\nonumber
\eea 
The following three statements are true :

Between, any zero of $P_j(x)$, and any element of 
B which is larger than this zero, there lies at
 least one zero of $P_{n\kp+j}(x)$.

Between any zero of $Q_j(x)$, and any element of A 
which is larger than this zero there lies at 
least one zero of $P_{n\kp+j}(x)$.

If $P_j(x)$ has a zero at any of the points in $B$ 
then $P_{n\kp+j}(x)$ has a zero at the same point.

\end{theorem}

{\it Proof:} 
From (4.7), $x_{i,n\kp+j}$ solves the following equation,
\bea
P_{n\kp}(x)P_j(x) = -
\frac{1}{\psi^2(x)}Q_{n\kp}(x)Q_j(x).
\nonumber
\eea
Using (4.5) and (4.6), we see that this is equivalent to,
\bea
\Delta_{\kp}T_n\left(\frac{P_\kp(x)}   
{\Delta_\kp}\right)P_j(x) =
U_n\left(\frac{P_\kp(x)}  
{\Delta_\kp}\right)Q_j(x)(1-x^2)  
\prod_{i=1}^{\kp-1}(x-\bt_i). 
\eea
For clarity, we illustrate a proof of the case of 
$\kp=3$ and $n=2$. The zeroes of $C(x)$ and 
$D(x)$ where,
\bea
C(x):=\Delta_{\kp}T_n
 \left(\frac{P_\kp(x)}{\Delta_\kp}\right)
\nonumber
\eea
and,
\bea 
D(x):=(1-x^2)U_n\left(\frac{P_\kp(x)}{\Delta_\kp}\right)
\prod_{i=1}^{\kp-1}(x-\bt_i)\nonumber
\eea
are indicated in figure 7. The crosses denote the 
zeros of $C(x)$ and the dots denote the zeros 
of $D(x)$. We see that they interlace. 
\begin{figure}[htb]
\centering
\includegraphics{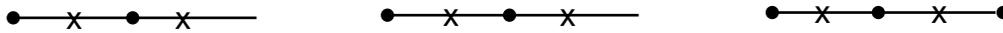}
\caption{{\it The zeroes of $C(x)$ and $D(x)$ for $\kp=3$ and $n=2$. The three line segments are 
the disjoint intervals that make up $E$. The crosses 
denote the zeroes of $C(x)$ and the dots, 
the zeroes of $D(x)$.}}
\end{figure}
In order to plot the remaining zeros of (4.24) we must 
add $j$ crosses and $j-1$ dots, 
corresponding to the zeros of $P_j(x)$ and $Q_j(x)$ 
respectively. We can distribute these without
 any restriction other than they must also interlace 
with each other. The result is a sequence of 
crosses and dots that alternate except for $j$ 
occurrences of two adjacent crosses and $j-1$ 
occurrences of two adjacent dots. Where we see 
two adjacent crosses, one of them is a zero of 
$P_j(x)$. Similarly with two dots, one of them is a 
zero of $Q_j(x)$. As an example we plot where
 the zeroes of both sides of (4.24) might lie, for 
$\kp=3$, $n=2$ and $j=2$ in figure 8. 
\begin{figure}[htb]
\centering
\includegraphics{dotcross2.eps}
\caption{}
\end{figure}

In order to graph the polynomials that pass through 
the zeros we need an initial condition. The 
condition is that both sides of (4.24) have the same 
sign at $x=-1+\epsilon$ where $\epsilon$ is 
a small positive constant. This is true because close 
to $x=-1$, $P_j(x)$ and $Q_j(x)$ must have 
opposite sign. Similarly, $T_n\left(\frac{P_\kp(x)}
{\Delta_\kp}\right)$ and 
$U_n\left(\frac{P_\kp(x)}{\Delta_\kp}\right)$ have 
the same sign if $\kp$ is even and have 
opposite sign if $\kp$ is odd. The product 
$\prod_{j=1}^{\kp-1}(x-\bt_i)$, for $x=-1+\epsilon$ 
is negative if $\kp$ is even and positive if $\kp$ is odd. 
$\Delta_\kp$ is a positive constant, 
and $(1-x^2)$ is also positive. Therefore both sides of 
(4.24) have the same sign at 
$x=-1+\epsilon$.

Therefore if we draw lines connecting the dots and the 
crosses that represent the tonicity of the
polynomials we will see the pattern shown in figure 9.
\begin{figure}[htb]
\centering
\includegraphics{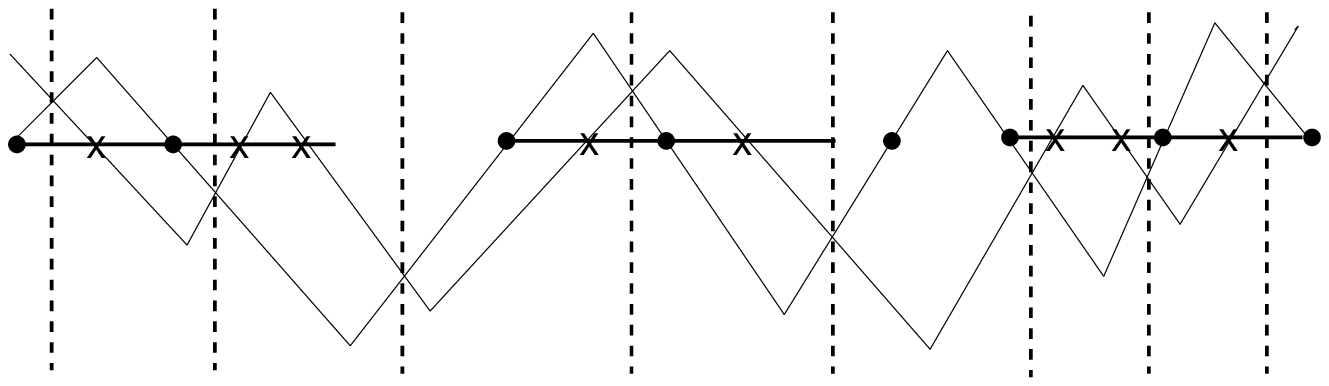}
\caption{The intersections of the two lines represent 
qualitatively where the zeroes of 
$P_{n\kp+j}(x)$ will lie relative to the zeroes of both 
sides of (4.24)}
\end{figure}

We are interested in where these polynomials intersect, 
as these points represent the zeros of 
$P_{n\kp+j}(x)$. Starting from the point $x=-1$ and 
moving along the interval $[-1,1]$, the 
following behaviour is observed. Before 
the first occurrence of two adjacent crosses the 
polynomials intersect between every dot-cross pair. After 
the appearance of the two adjacent 
crosses, the polynomials intersect between every 
cross-dot pair. This continues until we reach 
the first occurrence of two adjacent dots and then the 
intersection of the polynomials appears 
again in between every dot-cross pair. So every 
appearance of two crosses or two dots swaps the 
position of the intersection points between dot-cross 
pairs and cross-dot pairs. Since the 
appearance of two dots and two crosses alternate, 
after every appearance of two crosses there 
must be a point of intersection before the next dot. 
This corresponds to the statement in the 
theorem that between every zero of $P_j(x)$ and the 
next element from the set $B$ lies a zero of 
$P_{n\kp+j}(x)$. Similarly, after the appearance of 
two dots there must be an intersection point 
before the next cross. This is equivalent to saying 
that between every zero of $Q_j(x)$ and the 
next element from set $A$ also lies a zero of 
$P_{n\kp+j}(x)$.
$\Box$

{\bf Corollary}
{\it If $P_j(x)$ has a zero in $\bar{E_i}$ then so 
does $P_{n\kp+j}(x)$. ($\bar{E_i} := (\al_i,\bt_i))$}

{\it Proof:}
Referring to fig 9 we see that there is always a dot 
at the beta points. If there is a zero of $P_j(x)$ in
$\bar{E}_i$ then the previous theorem tells us that there
must be a zero of $P_{n\kp+j}$ between this zero and the
next $\bt$ point. Since $\bar{E}_i$ ends at the next $\bt$
point the zero of $P_{n\kp+j}$ is also contained in
$\bar{E}_i$. 
$\Box$

\newpage
\setcounter{section}{5}
\setcounter{equation}{0}
\setcounter{thm}{0}

{\bf Appendix A}

{\bf Genus one recurrence coefficients}

Here we show one method of determining the 
recurrence coefficients for the genus 1 case explicitly in
terms of the branch points $\al$ and $\bt$. It involves
manipulating the theta function expressions given in (1.8) and (1.9).
We see explicitly how the Jacobian elliptic functions define
the form of the coefficients.  From \cite{Chen-Lawrence}, for $g=1$ we have,
\bea
b_n = \frac{\bt-\al}{2}+\frac{1}{A}
\left[2\frac{\vt_3^{\pri}(u^+)} 
{\vt_3(u^+)}+\frac{\vt_3^{\pri}
((2n-3)u^+)}{\vt_3((2n-3)u^+)}
- \frac{\vt_3^{\pri}((2n-1)u^+)}
{\vt_3((2n-1)u^+)}\right], 
\nonumber
\eea
where
\bea
A = \frac{4K}{\sqrt{(1-\al)(1+\bt)}}\nonumber
\eea
and,
\bea
u^+ = \frac{1}{2K}\int_0^{\sqrt{\frac{\bt+1}{2}}}
\frac{dt}{{\sqrt {(1-t^2)(1-k^2t^2)}}}
\nonumber
\eea
with 
\bea
k^2 = \frac{2(\bt-\al)}{(1-\al)(1+\bt)}.\nonumber
\eea
Using,
\bea
\vt_3(v) = \vt_0(v\pm \frac{1}{2})\nonumber
\eea
and,
\bea
\theta_j(v) = \vt_j\left(\frac{v}{2K}\right)\nonumber
\eea
we write,
\bea
b_n = \frac{\bt-\al}{2}+\frac{2K}{A}\left[2\frac{\theta_0^{\pri}(2Ku^+ +K)}{\theta_0(2Ku^+ +K)}+\frac{\theta_0^{\pri}((2n-3)2Ku^+ -K)}{\theta_0((2n-3)2Ku^+ -K)}-\frac{\theta_0^{\pri}((2n-1)2Ku^+ +K)}{\theta_0((2n-1)2Ku^+ +K)}\right].\nonumber
\eea
The Jacobi function is defined to be
\bea
Z(w) = \frac{\theta_0^{\pri}(w)}{\theta_0(w)},\nonumber
\eea
and it satisfies the following formula,
\bea
Z(u+v)-Z(u-v)-2Z(v) = 
\frac{-2k^2\sn^2(u) \sn(v) \cn(v) \dn(v)}
 {1-k^2\sn^2u\sn^2v}.
\nonumber
\eea
If we make the choice,
\bea
u&=&(2n-2)2Ku^+ \nonumber\\
v&=&2Ku^+ +K,\nonumber
\eea
then we have,
\bea
b_n &=& \frac{2(\bt-\al)}
{\sqrt{(1-\al)(1+\bt)}}
\left[\frac{\sn^2((2n-2)2Ku^ +) \sn(2Ku^+ + K)
\cn(2Ku^+ +K)\dn(2Ku^+ +K)}
{1-k^2\sn^2((2n-2)2Ku^+)\sn^2(2Ku^+ +K)}\right]
\nonumber\\
&+& \frac{\bt-\al}{2}. 
\nonumber
\eea
Applying some basic identities $(k^{\pri}=\sqrt{1-k^2})$,
\bea
\sn(w+K) &=& \frac{\cn(w)}{\dn(w)}\nonumber\\
\cn(w+K) &=& -k^{\pri}\frac{\sn(w)}{\dn(w)}\nonumber\\
\dn(w+K) &=& \frac{-k^{\pri}}{\dn(w)}\nonumber
\eea
and,
\bea
\sn^2w+\cn^2w &=& 1\nonumber\\
k^2\sn^2w+\dn^2w &=& 1\nonumber
\eea
we can write,
\bea
b_n = (\bt-\al)\left[\frac{1}{2}
- \left(\frac{1+\al}{(\bt-\al)-(1+\bt)/\sn^2((2n-2)2Ku^+)}
\right)\right]
\nonumber
\eea
since,
\bea
\sn(2Ku^+) = \sqrt{\frac{\bt+1}{2}}.\nonumber
\eea
We have the following addition formula for 
$\sn(u)$ valid for any $u$ and $v$,
\bea
\sn(u+v) = \frac{\sn(u)\cn(v)\dn(v)
+ \sn(v)\cn(u)\dn(u)}{1-k^2\sn^2(u)\sn^2(v)}.
\nonumber
\eea
We can iterate this to obtain $\sn((2n-2)2Ku^+)$ 
in terms of $k$ and $\sn(2Ku^+)$. This gives us 
in principle, a method for evaluating the recurrence 
coefficients in terms of $\al$ and $\bt$. 
Owing to the nature of the addition formula, after 
just a few iterations the expressions for the 
coefficients become very large.

In \cite{Chen-Lawrence}, the $a_n$ are given as,
\bea
a_n = \left[\frac{\vt^{\pri}_1(0)}{A\vt_1(2u^+)}\right]^2
\frac{\vt_3((2n+1)u^+)\vt_3((2n-3)u^+)}
{\vt_3^2((2n-1)u^+)}\qquad n\geq 2. 
\nonumber
\eea
Using the following identity,
\bea
1 - k^2\sn^2(u)\sn^2(v) = \theta_0(0)^2
\frac{\theta_0(u+v)\theta_0(u-v)}
{\theta_0^2(u)\theta_0^2(v)}, 
\nonumber
\eea
with
\bea
u &=& (2n-1)2Ku^+ +K\nonumber\\
v &=& 4Ku^+ \nonumber
\eea
we have,
\bea
a_n &=& c\left(1-k^2\sn^2((2n-1)2Ku^+ +K)
\sn^2(4Ku^+)\right) \qquad n\geq 2. 
\nonumber\\
&=&c
\left(1-\frac{8(\bt-\al)} {(2 - \al + \bt)^2}
\sn^2((2n-1)2Ku^+ +K)\right). 
\nonumber
\eea

Here $c$ is a constant independent of n,
\bea
c=\left(\frac{\vt_1^{\pri}(0)\vt_0(2u^+)}
{A\vt_0(0)\vt_1(2u^+)}\right)^2.\nonumber
\eea
To determine $c$ explicitly, we evaluate 
$a_2$ by comparing coefficients of (1.13) for $n=2$ and 
then compare this with $a_2$ given above. 
In which case we find,
\bea
c = \frac{1}{16}(2-\al+\bt)^2.\nonumber
\eea
This value for $c$ was also determined in \cite{Chen-Lawrence} pg. 4695.
Note that in the higher genus cases, 
the non linear difference equations derived in section 2 can
be used to evaluate the recurrence coefficients.

\newpage
\setcounter{section}{6}
\setcounter{equation}{0}
\setcounter{thm}{0}

{\bf Appendix B}

{\bf g=2 Auxiliary polynomials}

Using the method outlined in section $2$ for $g=2$ we find,
\bea
\frac{G_3(x,n)}{h_{n-1}} &=& 
x^3-\frac{1}{2}(\al_1+\al_2
+ \bt_1+\bt_2)x^2\nonumber\\
&& \;\; -\Bigg(\frac{1}{2}
 + \frac{(\al_1^2+\al_2^2 + \bt_1^2 + \bt_2^2)}{8}
\nonumber\\
&& \qquad -\frac{1}{4}(\bt_1\bt_2 + \al_2\bt_1
+\al_2\bt_2+\al_1\al_2+\al_1\bt_1
+ \al_1\bt_2)-2a_n\Bigg)x
\nonumber\\
&&\;\;\;
+\;\frac{1}{16}
\Big(-(\al_1^3+\al_2^3
 + \bt_1^3+\bt_2^3) + 4(\al_1+\al_2+\bt_1 + \bt_2)  
 \nonumber\\
&&\;\;\;
+\;\bt_1\bt_2(\bt_1+\bt_2) + \al_2\bt_1(\al_2 + \bt_1)
+\al_2\bt_2(\al_2+\bt_2) + \al_1\al_2(\al_1+\al_2)   
\nonumber\\ &&\;\;\;
+\;\al_1\bt_1(\al_1 + \bt_1) + \al_1\bt_2(\al_1 + \bt_2)
\nonumber\\
&&\;\;\;
-\;2(\al_2\bt_1\bt_2+\al_1\bt_1\bt_2+\al_1\al_2\bt_1
+ \al_1\al_2\bt_2)
\nonumber\\
&&\;\;\;
-\;16a_n(\al_1 + \al_2 + \bt_1 + \bt_2-2(b_n +
b_{n+1}))\Big),   \qquad n \geq 2   
\nonumber
\eea
and,
\bea
\frac{S_2(x,n)}{h_n} &=& 
2x^2-\left[(\al_1+\al_2+\bt_1+\bt_2)-2b_{n+1}\right]x
\nonumber\\
&-&\frac{1}{4}
\Big(4 + \al_1^2 + \bt_1^2 + \al_2^2 + \bt_2^2-2(\al_1\al_2
+ \al_1\bt_1 + \al_2\bt_1 + 
\al_1\bt_2+\al_2\bt_2+\bt_1\bt_2)
\nonumber\\
&-&
8(a_{n+1} +a_n) +
4b_{n+1}(\al_1 + \al_2 + \bt_1 + \bt_2) - 8 b_{n+1}^2\Big),
\qquad n \geq 1.
\nonumber
\eea
\bea
G_3(x;1) = (x-b_1)(x-\al_1)(x-\al_2)\nonumber
\eea
and,
\bea
S_2(x;0) = (x-\al_1)(x-\al_2).\nonumber
\eea
Incidentally, this gives us the $\gamma_{1,2}(n)$ 
coefficients for the genus 2 case as,
\bea
\gamma_{1,2}(n) &=& \frac{1}{4}
(\al_1+\al_2+\bt_1+\bt_2)-2b_{n+1}
\nonumber\\
&\mp&\Bigg(8+3(\al_1^2+\al_2^2+\bt_1^2
+
\bt_2^2)-2(\al_1\al_2+\al_1\bt_1
+\al_2\bt_1+\al_1\bt_2+\al_2\bt_2+\bt_1\bt_2)
\nonumber\\
&-&16(a_n+a_{n+1})+4b_{n+1}(\al_1+\al_2+\bt_1
+\bt_2)-12b_{n+1}^2\Bigg)^{\frac{1}{2}},
\qquad n \geq 1.
\nonumber
\eea

\end{document}